\documentclass[amsfonts]{article}
\usepackage{amssymb}

\newtheorem{theorem}{Theorem}[section]
\newtheorem{corollary}[theorem]{Corollary}
\newtheorem{example}[theorem]{Example}

\newtheorem{lemma}[theorem]{Lemma}
\newtheorem{definition}[theorem]{Definition}
\newtheorem{remark}[theorem]{Remark}

\def\bR{\mathbb{R}}

\def\bC{\mathbb{C}}

\def\cB{\mathcal{B}}

\def\cD{\mathcal{D}}
\def\cE{\mathcal{E}}
\def\cF{\mathcal{F}}

\def\cH{\mathcal{H}}
\def\cI{\mathcal{I}}

\def\cM{\mathcal{M}}

\def\cR{\mathcal{R}}
\def\cS{\mathcal{S}}

\def\cL{\mathcal{L}}

\def\tM{\widetilde{\cM}}
\def\tmu{\widetilde{\mu}}

\begin{document}

\title{Linear SPDEs with harmonizable noise}

\author{Raluca M. Balan\footnote{Department of Mathematics and Statistics, University of Ottawa,
585 King Edward Avenue, Ottawa, ON, K1N 6N5, Canada. E-mail
address: rbalan@uottawa.ca} \footnote{Research supported by a
grant from the Natural Sciences and Engineering Research Council
of Canada.}}

\date{August 10, 2011}

\maketitle

\begin{abstract}
\noindent Using tools from the theory of random fields with stationary increments, we introduce a new class of processes which can be used as a model for the noise perturbing an SPDE. This type of noise (called harmonizable) is not necessarily Gaussian, but
it includes the spatially homogeneous Gaussian noise introduced in \cite{dalang99}, and the fractional noise considered in \cite{BT10}. We derive some general conditions for the existence of a random field solution of a linear SPDE with harmonizable noise, under some mild conditions imposed on the Green function of the differential operator which appears in this equation. This methodology is applied to the study of the heat and wave equations
(possibly replacing the Laplacian by one of its fractional powers), extending in this manner the results of \cite{BT10} to the case $H<1/2$.

\end{abstract}



\vspace{5mm}


\section{Introduction}

The theory of random processes, or random fields, with stationary increments was developed in the 1950's, the landmark references being \cite{doob53} and \cite{yaglom57}.
These processes have many interesting properties, most of which can be derived from their spectral representation as an integral with respect to a complex random measure with orthogonal increments.
In particular, various classes of self-similar processes with stationary increments, also called ``$H$-sssi processes'' ($H$ being the index of self-similarity), received special attention. These include the Hermite processes, which arise as limits in the non-central limit theorems (see \cite{dobrushin-major79}) and the generalized grey
Brownian motions (see \cite{mura-mainardi09}). The concept of operator scaling random field (OSRF) has been recently introduced as a replacement for self-similarity in the case of anisotropic fields. An explicit form of the spectral density of a Gaussian OSRF with stationary increments was obtained in \cite{clausel-vesel11}. An extension to the multivariate case was studied in \cite{li-xiao11}.

The {\em fractional Brownian motion} (fBm) is the most commonly known 
process with stationary increments.
There are two generalizations of the fBm in higher dimensions: the {\em fractional Brownian field} (fBf), and the {\em fractional Brownian sheet} (fBs).
The fBf has stationary increments, being the only isotropic Gaussian $H$-sssi process (see \cite{lindstrom93}, \cite{bonami-estrade03}).
The fBs was introduced in \cite{kamont96}, and was studied by many authors (see \cite{leger-pontier99}, \cite{bardina-florit05}). It does not have stationary increments, but can be studied using spectral analysis due to its representation as in integral with respect to a complex Gaussian measure. In the recent years, the fBs has been used increasingly as a model for the noise in stochastic analysis (see \cite{hu01}, \cite{hu-nualart-song11}).

The goal of this article is to introduce the tools necessary for the
study of stochastic partial differential equations (SPDEs) driven by a noise which has a harmonic structure similar to that of a random field with stationary increments, in both space and time. For this, an important step is to define a stochastic integral with respect to the noise.

In dimension $d=1$,
the
definition of a stochastic (or Wiener) integral with respect to
a process with stationary increments is usually taken for granted in the literature. In
the case of the fBm, this problem was studied thoroughly in
\cite{PT00} and \cite{PT01}, but these references did not explore
the fact that the spaces of Wiener integrands may contain
distributions. 
Implicitly used by several authors,
this fact has been recently proved in
\cite{jolis07} (in the case of the fBm) and
\cite{jolis10} (in the case of an arbitrary process with stationary increments).

By analogy with the case $d=1$, one may think that a natural space of Wiener integrands with respect to a random field $\{X_t\}_{t \in \bR^{d}}$ with stationary increments is the Hilbert space $\cH$ defined as the completion of the set $\cS$ of linear combinations of indicator functions with respect to the inner product $\langle 1_{[0,t]}, 1_{[0,s]}\rangle_{\cH}=E(X_t X_s)$.
However, in dimensions $d \geq 2$, even in the presence of the strong structure induced by the stationary increments, a collection of random variables indexed by the points in $\bR^d$ is not rich enough to produce a Hilbert space $\cH$ of integrands suitable for the study of a general SPDE.

Therefore, in the present article we assume that the space-time noise is a collection
$X=\{X(\varphi); \varphi \in \cD_{\bC}(\bR^{d+1})\}$ of random variables (indexed by the class $\cD_{\bC}(\bR^{d+1})$ of complex-valued infinitely differentiable functions on $\bR^{d+1}$ with compact support),
which can be represented as stochastic integrals with respect to a complex random measure $\tM$ with orthogonal increments on $\bR^{d+1}$.

In particular, the covariance of the noise is given by
$$E[X(\varphi)X(\psi)]=\int_{\bR^{d+1}} \cF \varphi(\tau,\xi) \overline{\cF \psi(\tau,\xi)} \Pi(d\tau,d\xi),$$ where $\cF \varphi$ denotes the Fourier transform of $\varphi$, and $\Pi$ is a symmetric tempered measure on $\bR^{d+1}$, which is related to the control measure of $\tM$.
This definition produces a rich enough space $\cH$, and induces a structure which is amenable to spectral analysis, being very similar to the structure of a random field with stationary increments on $\bR^{d+1}$. Moreover, this framework includes the Gaussian noise considered by other references (e.g. \cite{B11}, \cite{BT08}, \cite{BT10}, \cite{dalang99}, \cite{dalang-mueller03}, \cite{dalang-sanzsole05}), without being restricted to the Gaussian case. In particular, it provides a unifying framework for the study of equations with fractional noise in time, which covers simultaneously the cases $H>1/2$ and $H<1/2$.

An interesting case arises when $\Pi=\nu\times \mu$, where $\mu$ is a 
tempered measure on $\bR^d$, and 
$\tilde{\nu}(d\tau)=\tau^{-2}\nu(d\tau)$ is the spectral measure of a process $(Z_t)_{t \in \bR}$ with stationary increments and covariance function $R$. In this case, the noise can be identified with a collection $\{X_t(\varphi);t>0,\varphi \in \cD_{\bC}(\bR^d)\}$ of random variables which admit the representation mentioned above and have covariance:
$$E[X_t(\varphi)X_s(\psi)]=R(t,s) J(\varphi,\psi),$$ where
$J(\varphi,\psi)=\int_{\bR^d}\cF \varphi (\xi)\overline{\cF \psi(\xi)}\mu(d\xi)$.

In the Gaussian case, the process $\{X_t(\varphi);t \geq 0, \varphi \in \cD(\bR^d)\}$ was considered as a model for the noise driving an SPDE by many references. These references considered {\em only the case} when the Fourier transform of $\mu$ in $\cS'(\bR^d)$ is a locally integrable function $f$ (or a measure $\Gamma$), and hence $J$ can be represented as:
\begin{equation}
\label{rep-space-covariance1} J(\varphi,\psi)=\int_{\bR^d}
\int_{\bR^d}\varphi(x)\psi(y)f(x-y)dxdy, \quad \forall
\varphi,\psi \in {\cS}(\bR^d).
\end{equation}
See e.g. \cite{dalang99}, \cite{dalang-mueller03}, \cite{dalang-sanzsole05} for the case when $Z$ is a Brownian motion, and \cite{BT08}, \cite{BT10}, \cite{B11} for the case when $Z$ is a fBm of index $H>1/2$. In general, the Fourier transform of $\mu$ in $\cS'(\bR^d)$  may not be a locally integrable function (or a measure). This is the case for instance, when $d=1$ and $\mu(d\xi)=|\xi|^{-\alpha}$ with $\alpha \in (-1,0)$.
The fact that representation (\ref{rep-space-covariance1}) is not needed in the present article leads to a more general class of noise processes, even in the Gaussian case.

This article is organized as follows.
In Section \ref{background-section} we review some background material related to random fields with stationary increments.
In Section \ref{noise-section}, we introduce the ``harmonizable'' noise $X$ (described above), and we identify several classes of deterministic integrands with respect to this noise.
In Section \ref{equations-section}, we define the random-field solution of a linear SPDE $Lu=\dot{X}$ driven by the harmonizable noise $X$, and we show that, under some mild regularity conditions imposed on the fundamental solution $G$ of the equation $Lu=0$, the necessary and sufficient condition for the existence of this solution is:
\begin{equation}
\label{gen-cond}
\int_{\bR^{d+1}}\left|\int_0^t e^{-i \tau s} \cF G(s,\cdot)(\xi)ds \right|^2 \Pi(d\tau,d\xi)<\infty,
\end{equation}
where $\cF G(s,\cdot)$ denotes the Fourier transform of the distribution $G(s,\cdot) \in \cS'(\bR^d)$.

Section \ref{parabolic-section} and Section \ref{hyperbolic-section} are dedicated, respectively, to the parabolic case $L=\frac{\partial}{\partial t}-\cL$, and the hyperbolic case $L=\frac{\partial^2}{\partial t^2}-\cL$, where $\cL=-(-\Delta)^{\beta/2}$ for some $\beta>0$, assuming that $\Pi=\nu \times \mu$ (with measures $\nu$ and $\mu$ as above). In the parabolic case,
we show that condition (\ref{gen-cond}) is equivalent to:
\begin{equation}
\label{parabolic-cond} \int_{\bR^d}
\int_{\bR}\frac{1}{\tau^2+[1+\Psi(\xi)]^2}\nu(d\tau)\mu(d\xi)<\infty,
\end{equation}
where $\Psi(\xi)=|\xi|^{\beta}$. In the hyperbolic case, condition (\ref{gen-cond}) is equivalent to:
\begin{equation}
\label{hyperbolic-cond2} \int_{\bR^d}\frac{1}{\sqrt{1+\Psi(\xi)}}
\int_{\bR} \frac{1}{\tau^2+1+\Psi(\xi)}\nu(d\tau)
\mu(d\xi)<\infty.
\end{equation}

When $(Z_t)_{t \in \bR}$ is the fBm of index $H \in (0,1)$, conditions (\ref{parabolic-cond}) and (\ref{hyperbolic-cond2}) become $$\int_{\bR^d}\left(\frac{1}{1+\Psi(\xi)}\right)^{2H}\mu(d\xi)<\infty, \ \mbox{respectively} \ \int_{\bR^d}\left(\frac{1}{1+\Psi(\xi)}\right)^{H+1/2}\mu(d\xi)<\infty,$$
extending the results of \cite{BT10} to the case $H<1/2$.

The appendix collects some auxiliary results, which are invoked in the article.
In Appendix A, we review the construction of the stochastic integral with respect to a complex random measure with orthogonal increments, and list some of its properties. Appendix B contains some elementary estimates related to the fBm case. Appendix C shows that some technical conditions (which are imposed for treating the hyperbolic equation) are verified when the measure $\nu$ is the Fourier transform of the Riesz kernel or the Bessel kernel.

We conclude the introduction with a few words about the notation. We let $\cD_{\bC}(\bR^{d})$
be the set of complex-valued infinitely differentiable functions on $\bR^{d}$ with compact support.
For any $p>0$, we denote by $L_{\bC}^p(\bR^d,\mu)$ the space of complex-valued
functions $\varphi$ on $\bR^d$ such that $|\varphi|^p$ is integrable with respect to the measure $\mu$. When $\mu$ is the Lebesgue measure, we simply write $L_{\bC}^p(\bR^d)$.
We denote by $L_{\bC}^p(\Omega)$ the set of complex-valued random variables $X$, defined on a probability space $(\Omega,\cF,P)$, such that $E|X|^p<\infty$.
We let ${\cS}_{\bC}(\bR^d)$ be the set of complex-valued rapidly decreasing
infinitely differentiable functions on $\bR^d$.
Similar notations are used for spaces of real-valued elements, with the subscript $\bC$ omitted.

We denote by $\cF \varphi(\xi)=\int_{\bR^d}e^{-i \xi \cdot x} \varphi(x)dx$ the
Fourier transform of a function $\varphi \in L^1(\bR^d)$. (We
use the same notation $\cF$ for the Fourier transform of functions on $\bR,\bR^d$ or $\bR^{d+1}$. Whenever there is a risk of confusion, the notation will be clearly
specified.) We denote by $x \cdot y=\sum_{i=1}^{d}x_iy_i$ the inner
product in $\bR^d$ and by $|x|=(x \cdot x)^{1/2}$ the Euclidean norm
in $\bR^d$. 
A bounded rectangle in $\bR^d$ is a set of the form $(x,y]=\{z \in \bR^d; x_i < z_i \leq y_i \ \mbox{for all} \ i=1, \ldots,d\}$. We denote by $\cR_d$ the class of bounded rectangles in $\bR^d$ and by $\cB_{b}(\bR^d)$ the class of bounded Borel sets in $\bR^{d}$.

\section{Background}
\label{background-section}


In this section, we review the basic elements of the theory of random processes and random fields with stationary increments. We refer the reader to Sections 23 and 25.2 of \cite{yaglom87} for more details. 

\begin{definition}
{\rm A real-valued random field $X=(X_t)_{t \in
\bR^d}$, defined on
a probability space $(\Omega,\cF,P)$, has {\em (wide-sense) stationary increments},
if $E(X_t^2)<\infty$ for any $t \in \bR^d$, and the
increments of the process have homogenous first and second order
moments, i.e. $E(X_{t+s}-X_{s})=m(t)$ for any $t,s \in \bR^d$, and
\begin{equation}
\label{homog-cov} E|X_{t+s}-X_{s}|^2=V(t) \quad \mbox{for any} \
t,s \in \bR^d.
\end{equation}
}
\end{definition}
In particular, a stationary random field has stationary increments.

In the present article, we work only with zero-mean processes, i.e. we assume that $m(t)=0$ for all $t \in \bR^d$.
From (\ref{homog-cov}), it follows that for any $s,t,u \in \bR^d$,
$$E[(X_{u+t}-X_{u})(X_{u+s}-X_{u})] =
\frac{1}{2}[V(t)+V(s)-V(t-s)]=:R(t,s).$$
(To see this, write $(X_{u+t}-X_{u})(X_{u+s}-X_{u})=
\frac{1}{2}(|X_{u+t}-X_{u}|^2+|X_{u+s}-X_{u}|^2
-|X_{u+t}-X_{u+s}|^2)$.)
On the other hand, the covariance between arbitrary increments
$X_{u+t}-X_u$ and $X_{v+s}-X_{v}$ depends on $t,s$ and $u-v$:
$$E[(X_{u+t}-X_u)(X_{v+s}-X_{v})]=R(t,s+v-u)-R(t,v-u).$$
This follows by writing $X_{v+s}-X_v=(X_{v+s}-X_u)-(X_v-X_u)$.

The function $R$ is called the {\em covariance function} of $X$
and plays an important role in the analysis of the process $X$,
via the spectral theory.
More precisely, if $R$ is continuous, then the process $X$ and
its covariance function admit the following {\em spectral representations}:
(see (4.351) and (4.345) of \cite{yaglom87})
\begin{eqnarray}
\label{spectral-rep-stat-incr}X_t&=&X_0+t \cdot Y +\int_{\bR^d
\verb2\2 \{0\}} (e^{-i \tau \cdot t}-1) \tM(d\tau)\\
\nonumber
R(t,s)&=& \Sigma t \cdot s +\int_{\bR^d \verb2\2 \{0\}} (e^{-i\tau
\cdot t}-1)(e^{i\tau \cdot s}-1)\tmu(d\tau).
\end{eqnarray}
The elements which appear in these representations are the
following:
\begin{itemize}
\item $\tM$ denotes a collection $\{\tM(A); A \in \cR_d^0 \}$ of
zero-mean random variables in $L_{\bC}^2(\Omega)$, indexed by the class
$\cR_d^0$ of rectangles in $\bR^d$ which do not include $0$, which satisfy the following conditions:\\
{\em (i) (additivity)} $\tM(A \cup B)=\tM(A)+\tM(B)$ a.s. for disjoint sets $A,B \in \cR_d^0$; \\
{\em (ii) (orthogonality)} $E[\tM(A)\overline{\tM(B)}]=0$ for any disjoint sets $A,B \in \cR_d^0$; \\
{\em (iii) (symmetry)} $\tM(-A)=\overline{\tM(A)}$ a.s. for any $A \in \cR_d^0$; \\
{\em (iv) (control measure)} $E|\tM(A)|^2=\tmu(A)$ for any $A \in \cR_d^0$.

\item $\tmu$ is a symmetric measure on $\bR^d \verb2\2\{0\}$ which
satisfies the condition:
\begin{equation}
\label{cond-tmu} \int_{\bR^d \verb2\2
\{0\}}\frac{|\tau|^2}{1+|\tau|^2}\tmu(d\tau)<\infty.
\end{equation}
The measure $\tmu$ is called the {\em spectral measure} of $X$.

\item $Y$ is a real-valued centered $d$-dimensional random vector
with covariance matrix $\Sigma$, which is orthogonal to $\tM(A)$ for all $A \in \cR_d^0$, i.e. $E[Y \overline{\tM(A)}]=0$ for all $A \in \cR_d^0$. We may define $\tM(\{0\})=Y$ and $\tmu(\{0\})=E|Y|^2$.
\end{itemize}

Condition (\ref{cond-tmu}) ensures that the function $\varphi(\tau)=e^{-i \tau \cdot t}-1$ is in $L_{\bC}^2(\bR^d,\tmu)$, and hence the integral of (\ref{spectral-rep-stat-incr}) is well-defined.

In the present article, we assume that $X_0=0$ and $Y=0$.
We consider an object $\cM$ as above, except that its indexing class is extended to all rectangles (which may include $0$), and the measure $\tmu$ is not required to satisfy (\ref{cond-tmu}).
More precisely, we introduce the following definition. 

\begin{definition}
\label{def-tM}
{\rm Let $\tmu$ be a symmetric measure on $\bR^d$. A collection $\tM=\{\tM(A); A \in \cR_d \}$ of zero-mean random variables in $L_{\bC}^2(\Omega)$ which satisfy conditions {\em (i)-(iv)} above for any sets $A,B \in \cR_d$ is called a {\em complex random measure on $\bR^d$ with orthogonal increments and control measure $\tmu$}.}
\end{definition}

The stochastic integral with respect to $\tM$ is defined as the mean-square limit of a sequence of Riemann-Stieltjes integrals. More precisely, for any function $\varphi \in L_{\bC}^2(\bR^d,\tilde{\mu})$ which is continuous on compact sets,
we define the following $\bC$-valued random variable:
\begin{eqnarray*}
\int_{\bR^d} \varphi(\tau)\tM(d\tau)
&=&\lim_{n \to \infty}\{\lim_{\max_j|A_j| \to 0} \sum_{j=1}^{k_n}\varphi(\tau_j)\tM(A_j) \} \quad \mbox{in} \ L_{\bC}^2(\Omega),
\end{eqnarray*}
where $(A_j)_{1 \leq j \leq k_n} \subset \cR_d$ is a partition of $[-N_n,N_n]$ with $N_n=(n, \ldots,n) \in \bR^d$, $|A_j|$ is the length of the longest edge of $A_j$, and $\tau_j \in A_j$ is  arbitrary. For the sake of completeness, the details of this construction are given in Appendix A. 


\vspace{1mm}

We now examine the case $d=1$.
Let $Z=(Z_t)_{t \in \bR}$ be a real-valued
centered random process with stationary increments and spectral representation (\ref{spectral-rep-stat-incr}) with $Y=0$ and spectral measure
$\tilde{\nu}$. We assume that $Z_0=0$.

Let
$\nu(d\tau)=\tau^2 \tilde{\nu}(d\tau)$. Note that $\nu$ is a symmetric measure on $\bR$ which satisfies the condition:
\begin{equation}
\label{nu-cond} K:=\int_{\bR}\frac{1}{1+\tau^2}\nu(d\tau)<\infty.
\end{equation}
The measure $\nu$ is sometimes also called the spectral measure of $Z$.

Since $\cF 1_{[0,t]}(\tau)=\int_0^t e^{-i \tau s}ds=(e^{-i\tau
t}-1)/(-i \tau)$, it follows that:
\begin{equation}
\label{time-covariance} R(t,s)=E(Z_t Z_s)=\int_{\bR} \cF 1_{[0,t]}(\tau)
\overline{\cF 1_{[0,s]}(\tau)} \nu(d\tau).
\end{equation}
Note that relation (\ref{time-covariance}) does not hold in dimensions $d \geq 2$.



Let $\cS$ be the set of elementary functions on $\bR$, i.e. real linear
combinations of indicator functions $1_{[0,t]}, t \in \bR$. We endow
$\cS$ with the inner product:
$$\langle 1_{[0,t]}, 1_{[0,s]}\rangle_{\Lambda}=R(t,s).$$

Let $\Lambda$ be the completion of $\cS$ with respect to $\langle
\cdot, \cdot \rangle_{\Lambda}$. The map $1_{[0,t]} \mapsto Z_t$
is an isometry between $\cS$ and $L^2(\Omega)$, which is
extended to $\Lambda$. This extension is denoted by $\varphi
\mapsto Z(\varphi)=\int_{\bR}\varphi(t)dZ_t$. The
space $\Lambda$ contains distributions. More precisely, if
$\varphi \in \cS'(\bR)$ is such that $\cF \varphi$ is a function
and $\int_{\bR}|\cF \varphi(\tau)|^2 \nu(d\tau)<\infty$, then
$\varphi\in \Lambda$ (see Theorem 3.2 of \cite{jolis10}).

\begin{remark}
\label{Lambda-rem} 
{\rm 
By Theorem 2.3 of \cite{jolis10}
and Proposition 2.5 of \cite{jolis10}, $\Lambda$ coincides with
the completion of $\cD(\bR)$ with respect to the inner
product:
$$\langle \varphi, \psi \rangle_{\Lambda}=\int_{\bR} \cF
\varphi(\tau)\overline{\cF \psi(\tau)}\nu(d\tau).$$
}
\end{remark}

In general, the Fourier transform in $\cS'(\bR)$ of the tempered
measure $\nu$ is a distribution (denoted by $\cF \nu$), which
coincides with $(1/2) V''$, $V''$ being the second order
distributional derivative of the continuous variance function $V$
(see Section 2.2 of \cite{jolis10}). In some cases, $\cF \nu$ is
a locally integrable function $\rho$. By the Fourier inversion theorem on $\cS(\bR)$, this is equivalent to saying that $\nu$ is the Fourier transform of $(2\pi)^{-1}\rho$ in $\cS'(\bR)$.
In this case, (see Lemma 5.6 of \cite{KX09})
\begin{equation}
\label{R-repr} R(t,s)=\int_{0}^{t} \int_{0}^{s}\rho(u-v)du dv.
\end{equation}
Relation (\ref{R-repr}) is useful in applications, but is not
needed in the present article.

When $\nu$ is finite,
$\rho(t)=\int_{\bR}e^{-i \tau t} \nu(d\tau)$, and a process $(Z_t)_{t \in \bR}$ with
stationary increments and covariance $R$ given by (\ref{R-repr}) can be represented as
 $Z_t=\int_0^t Y_s ds$, where $(Y_t)_{t \in \bR}$ is a
zero-mean stationary process with 
$E(Y_0 Y_t)=\rho(t)$.

\begin{example}
\label{Riesz-ex} {\rm Let $\nu(d\tau)=|\tau|^{-\gamma}$ for some
$\gamma \in (-1,1)$. 
Then $\nu$ is the Fourier transform of a locally integrable function if and only if $\gamma \in (0,1)$ (see Lemma 3.1 of \cite{jolis07}).
This function is the Riesz kernel $\rho(t)=c_{\gamma}|t|^{-(1-\gamma)}$
(see Lemma 1, Chapter V of
\cite{stein70}).
In the Gaussian case, the parametrization $\gamma=2H-1$ with $H \in (0,1)$ corresponds
to the case when $Z$ is a fBm of index $H$.
}
\end{example}

\begin{example}
\label{Bessel-ex} {\rm Let $\nu(d\tau)=(1+|\tau|^2)^{-\gamma/2}$
for some $\gamma >-1$. 
When $\gamma>0$, $\nu$
is the Fourier transform of the Bessel kernel (see Proposition 2, Chapter V of \cite{stein70})
$$\rho(t)=c_{\gamma}\int_{0}^{\infty} w^{(\gamma-1)/2-1}e^{-w} e^{-t^2/(4w)}dw.$$
The measure $\nu$ is finite if and only if $\gamma>1$. If $\gamma=2$, then
$\rho(t)=
c_1 \exp(-c_2|t|)$ for some constants $c_1,c_2>0$ (see Exercise
2.3.4, Chapter 8 of \cite{K02}), and a process $(Z_t)_{t \in \bR}$ with stationary
increments and spectral measure $\nu(d\tau)=\tau^{-2}(1+\tau^2)^{-1}d\tau$
can be represented as $Z_t=\int_0^t Y_s ds$, where $(Y_t)_{t \in
\bR}$ is a zero-mean stationary process with
$E(Y_0Y_t)=\rho(t)$. In the Gaussian case, $(Y_t)_{t
\in \bR}$ is a zero-mean Ornstein-Uhlenbeck process.

}
\end{example}

\section{The harmonizable noise}
\label{noise-section}

Let $\Pi$ be a (symmetric) tempered measure on $\bR^{d+1}$, i.e.
\begin{equation}
\label{cond-mu}\int_{\bR^{d+1}} \left( \frac{1}{1+\tau^2+|\xi|^2}\right)^k
\Pi(d\tau,d\xi)<\infty \quad \mbox{for some} \ k>0.
\end{equation}

Let $\tM=\{\tM(A); A \in \cR_{d+1}\}$ be a complex random measure on $\bR^{d+1}$, with orthogonal increments and control measure $\widetilde{\Pi}$ defined by:
$$\widetilde{\Pi}(d\tau,d\xi)=\frac{1}{\tau^2 \xi_1^2 \ldots \cdot \xi_d^2}\Pi(d\tau,d\xi) \quad \mbox{for} \ \tau \in \bR,\xi \in \bR^d.$$

For every $\varphi \in \cS_{\bC}(\bR^{d+1})$, we define a $\bC$-valued random variable by:
\begin{equation}
\label{def-X(varphi)}
X(\varphi)=\int_{\bR^{d+1}}\cF \varphi(\tau,\xi)(-i \tau)\prod_{j=1}^{d}(-i \xi_j) \tM(d\tau,d\xi).
\end{equation}

To see that $X(\varphi)$ is well-defined in $L_{\bC}^2(\Omega)$, we use the fact that $\cF \varphi\in \cS_{\bC}(\bR^{d+1})$. Hence
$\sup_{\tau \in \bR,\xi \in \bR^d}(1+\tau^2+|\xi|^2)^n|\cF \varphi(\tau,\xi)|^2<\infty$ for any $n>0$.
Taking $n=k$ and using (\ref{cond-mu}), it follows that
$$\int_{\bR^{d+1}}|\cF \varphi(\tau,\xi)|^2 \tau^2 \prod_{j=1}^{d}\xi_j^2 \ \tilde{\Pi}(d\tau,d\xi)=\int_{\bR^{d+1}}|\cF \varphi(\tau,\xi)|^2  \Pi(d\tau,d\xi)<\infty,$$ i.e. the function $\psi(\tau,\xi)=
\cF \varphi(\tau,\xi)(-i \tau)\prod_{j=1}^{d}(-i \xi_j)$ belongs to $L_{\bC}^2(\bR^{d+1},\tilde{\Pi})$. Since $\psi$ is also continuous on compacts, $X(\varphi)$ is well-defined (see Appendix \ref{app-stoch-integr}).

If $\varphi$ is real-valued, then $X(\varphi)$ is real-valued, due to the symmetry of $\tM$ (see relation (\ref{M-real-valued}), Appendix A).
Using (\ref{cov-M}), we see that for any $\varphi, \psi \in \cS_{\bC}(\bR^{d+1})$,
$$E[X(\varphi)X(\psi)]=\int_{\bR^{d+1}} \cF \varphi(\tau,\xi)\overline{\cF \psi(\tau,\xi)}\Pi(d\tau,d\xi)=:\cI(\varphi,\psi).$$

The functional $\cI$ is bilinear, symmetric, and non-negative definite.
To see this, note that for any $m \geq 1, a_1, \ldots,a_m \in \bC$ and $\varphi_1,
\ldots, \varphi_m \in \cS_{\bC}(\bR^{d+1})$,
$$\sum_{i,j=1}^{m}a_i \overline{a_j}\cI(\varphi_i,\varphi_j)=\int_{\bR^{d+1}}\left|\sum_{i=1}^{m}a_i \cF \varphi_i(\tau,\xi)\right|^2 \Pi(d\tau,d\xi) \geq 0.$$

We identify any functions $\varphi,\psi \in \cS_{\bC}(\bR^{d+1})$
for which $\cI(\varphi-\psi,\varphi-\psi)=0$, so that $\cI$ is an inner product.

\begin{remark}
{\rm If the Fourier transform of $\Pi$ in $\cS'(\bR^{d+1})$ is a locally integrable 
function $F:\bR^{d+1} \to [0,\infty)$,
then
$$\cI(\varphi,\psi)=\int_{\bR^{2}}
\int_{\bR^{2d}}\varphi(t,x)\psi(s,y)F(t-s,x-y)dxdydt ds,$$
for any $\varphi,\psi \in {\cS}(\bR^{d+1})$. This representation is not needed in the present work.
}
\end{remark}

In what follows, we work with functions $\varphi \in \cD_{\bC}(\bR^{d+1})$.
More precisely, we introduce the following definition.

\begin{definition}
{\rm The process $\{X(\varphi);\varphi \in \cD_{\bC}(\bR^{d+1})\}$ defined by (\ref{def-X(varphi)}) is called a {\em harmonizable noise} with driving measure $\tM$ and covariance measure $\Pi$.}
\end{definition}

We endow $\cD_{\bC}(\bR^{d+1})$ with the inner product $\langle \cdot,\cdot \rangle_{\cH}$ defined by:
$$\langle \varphi,\psi \rangle_{\cH}=\cI(\varphi,\psi).$$
and we let $\cH$ be the completion of $\cD_{\bC}(\bR^{d+1})$ with respect to this inner product.

The map $\varphi \mapsto X(\varphi)$ is an isometry from $\cD_{\bC}(\bR^{d+1})$ to $L_{\bC}^2(\Omega)$, which is extended to $\cH$.
The Hilbert space $\cH$ contain distributions.
By abuse of notation, we write $$X(\varphi)=\int_{\bR}
\int_{\bR^d}\varphi(t,x)X(dt,dx)$$ even if $\varphi$ is a
distribution. This stochastic integral is well-defined only if
$\varphi \in \cH$.

We will need an alternative representation for the inner product
in $\cH$. For any $\varphi_1, \varphi_2 \in \cD_{\bC}(\bR^{d+1})$, if we denote by
$$\phi_{\xi}^{(i)}(t)=\cF \varphi_i(t, \cdot)(\xi)$$ the
Fourier transform of the function $\varphi_i(t, \cdot)$, and by
$\cF\phi_{\xi}^{(i)}$ the Fourier transform of the function $t
\mapsto \phi^{(i)}_{\xi}(t)$ for $i=1,2$, then
\begin{equation}
\label{H-product-E} \langle \varphi_1,\varphi_2
\rangle_{\cH}=\int_{\bR^{d+1}} \cF\phi_{\xi}^{(1)}(\tau)
\overline{\cF\phi_{\xi}^{(2)}(\tau)} \Pi(d\tau,d\xi)=:\langle
\varphi_1,\varphi_2 \rangle_{0}.
\end{equation}

 The following result is a modified version of Lemma 3.7 of
\cite{BT08}.

\begin{lemma}
\label{modified-BT-lemma} Let $\mu$ be a tempered measure on
$\bR^{d}$. For every $A \in \cB_b(\bR^d)$, there exists
a sequence $(\psi_n)_{n \geq 1} \subset \cD_{\bC}(\bR^d)$ such
that
$$\int_{\bR^d}|1_{A}(\xi)-\cF \psi_n(\xi)|^2 \mu(d\xi) \to 0.$$
\end{lemma}

\noindent {\bf Proof:} As in the proof of Lemma 3.7 of
\cite{BT08}, let $\eta_n(x)=n^{d}\eta(nx)$, where $\eta \in
\cD_{\bC}(\bR^d)$ is such that $\eta \geq 0$ and $\int
\eta(x)dx=1$. Setting $\phi_n=1_{A}* \eta_n$, we note that $\phi_n
\in \cD_{\bC}(\bR^d)$, $\phi_n \to 1_{A}$ uniformly, and ${\rm
supp} \ \phi_n \subset K$ for all $n$, where $K$ is a compact set
in $\bR^d$. By the dominated convergence theorem,
$$\int_{\bR^d}|\phi_n(\xi)-1_{A}(\xi)|^2 \mu(d\xi) \to 0.$$

For the second part of the proof, we argue as on p. 613 of
\cite{jolis10}. Since the Fourier transform is a homeomorphism
from $\cS(\bR^d)$ onto itself, and $\cD_{\bC}(\bR^d)$ is dense
in $\cS(\bR^d)$, $\cF(\cD_{\bC}(\bR^d))$ is dense in
$\cS_{\bC}(\bR^d)$. Since $\mu$ is a tempered measure, the
convergence in $\cS_{\bC}(\bR^d)$ implies the convergence in
$L_{\bC}^2(\bR^d,\mu)$. Hence, there exists a sequence $(\psi_n)_{n \geq
1} \subset \cD_{\bC}(\bR^d)$ such that:
$$\int_{\bR^d}|\phi_n(\xi)-\cF \psi_n(\xi)|^2 \mu(d\xi) \to 0.$$
The conclusion follows. $\Box$

Next, we give a criterion for checking that $\varphi \in \cH$,
which is a generalization of Theorem 2.1 of \cite{BT10}.

\begin{theorem}
\label{varphi-in-H} Let $\bR \ni t \mapsto \varphi(t,\cdot) \in
\cS'(\bR^d)$ be a deterministic function such that
$\cF \varphi(t,\cdot)$ is a function for all $t \in
\bR$.
Suppose that:\\
(i) for all $\xi \in \bR^d$, the function $t \mapsto \cF
\varphi(t,\cdot)(\xi)=:\phi_{\xi}(t)$ belongs to $L_{\bC}^1(\bR)$; \\
(ii) the function $(t,\xi) \mapsto \cF \varphi(t,\cdot)(\xi)$ is
measurable on $\bR \times \bR^d$.

For any $\xi \in \bR^d$, we denote by $\cF \phi_{\xi}$ the Fourier
transform of $\phi_{\xi}$. If
\begin{equation}
\label{norm-varphi-0-finite}\|\varphi\|_{0}^2:=\int_{\bR^{d+1}}
|\cF \phi_{\xi}(\tau)|^2 \Pi(d\tau,d\xi)<\infty,
\end{equation}
then $\varphi \in \cH$ and $\|\varphi
\|_{\cH}^2=\|\varphi\|_{0}^{2}$. In particular, the stochastic
integral of $\varphi$ with respect to $X$ is well-defined and is given by:
\begin{equation}
\label{rep-X(varphi)}
X(\varphi)=\int_{\bR^{d+1}} \cF \phi_{\xi}(\tau) (-i \tau) \prod_{j=1}^{d}(-i \xi_j)\tM(d\tau,d\xi).
\end{equation}
\end{theorem}

\noindent {\bf Proof:} The argument is a simplified version of the proof of Theorem 2.1 of \cite{BT10}. Due to
(\ref{H-product-E}) and the definition of $\cH$, we have to
show that for any $\varepsilon>0$ there exists $l=l_{\varepsilon}
\in \cD_{\bC}(\bR^{d+1})$ such that
\begin{equation}
\label{norm-a-b}\|\varphi-l\|_{0}^{2}=\int_{\bR^{d+1}}|\cF
\phi_{\xi}(\tau)-\cF \psi_{\xi}(\tau)|^2\Pi(d\tau,d\xi)<\varepsilon^2,
\end{equation}
where $\cF \psi_{\xi}$ denotes the Fourier transform of the
function $t \mapsto \psi_{\xi}(t)=\cF l(t, \cdot)(\xi)$.

We denote $a(\tau,\xi)=\cF \phi_{\xi}(\tau)$. It will be shown
later that:
\begin{equation}
\label{a-is-measurable}\mbox{the function $(\tau,\xi) \mapsto
a(\tau,\xi)$ is measurable.}
\end{equation}

Due to (\ref{norm-varphi-0-finite}) and (\ref{a-is-measurable}),
it follows that $a \in L_{\bC}^2(\bR^{d+1},\Pi)$. By applying
Theorem 19.2 of \cite{billingsley95} to the real and imaginary parts of the function $a$,
we infer that there exists a simple function $h$ on $\bR^{d+1}$ such that
\begin{equation}
\label{norm-a-h} \int_{\bR^{d+1}}
|a(\tau,\xi)-h(\tau,\xi)|^2
\Pi(d\tau,d\xi)<\frac{\varepsilon^2}{4}.
\end{equation}

Here, a simple function is a complex linear combination of indicator
functions of the form $1_{B}$, where $B$ is a Borel subset of
$\bR^{d+1}$ which is not necessarily bounded. From the proof of Theorem 19.2 of
\cite{billingsley95}, we see that the Borel sets $B$ which appear
in the definition of the function $h$ 
satisfy
$\Pi(B)<\infty$. The measure $\Pi$
is locally finite, and hence $\sigma$-finite. By applying Theorem 11.4.(ii) of
\cite{billingsley95}, each of these Borel sets can be
approximated by a finite union of bounded rectangles in
$\bR^{d+1}$. Therefore, we may assume that each set $B$ is bounded.

By Lemma \ref{modified-BT-lemma}, there exists a function $l \in
\cD_{\bC}(\bR^{d+1})$ such that
\begin{equation}
\label{norm-b-h}\int_{\bR^{d+1}}|1_{B}(\tau,\xi)-\cF
l (\tau,\xi)|^2\Pi(d\tau,d\xi)<\frac{\varepsilon^2}{4}.
\end{equation}

Relation (\ref{norm-a-b}) follows from (\ref{norm-a-h}) and
(\ref{norm-b-h}), since $\cF l(\tau,\xi)=\cF \psi_{\xi}(\tau)$.

We now prove (\ref{a-is-measurable}). Since $(t,\xi) \mapsto
\phi_{\xi} (t)$ is measurable, by applying Theorem 13.5 of
\cite{billingsley95} to the real and imaginary part of
$\phi_{\xi}(t)$ we conclude that there exists a sequence $(l_n)_n$
of complex-valued simple functions such that $l_n(t,\xi) \to
\phi_{\xi}(t)$ and $|l_n(t,\xi)| \leq |\phi_{\xi}(t)|$ for any $n
\geq 1$. We let $a_{l_n}(\tau,\xi)=\cF l_n(\cdot,\xi)(\tau)$. By
the dominated convergence theorem,
$$a_{l_n}(\tau,\xi)-a(\tau,\xi)=\int_{0}^{\infty}e^{-i\tau t} (l_n(t,\xi)dt-\phi_{\xi}(t))dt \to 0.$$
If $l(t,\xi)=1_{(c,d]}(t)1_{A}(\xi)$ is an elementary function,
then $a_{l}(\tau,\xi)=\cF l(\cdot,\xi)(\tau)=\cF
1_{(c,d]}(\tau)1_{A}(\xi)$ is clearly measurable. Since
$a_{l_n}(\tau,\xi)$ is measurable for any $n$, it follows that
$a(\tau,\xi)$ is measurable.

To prove (\ref{rep-X(varphi)}), we denote by $Y$ the random variable on the right hand side of (\ref{rep-X(varphi)}), which is well-defined due to (\ref{norm-varphi-0-finite}) and hypothesis {\em (i)}. By definition, $X(\varphi)=\lim_{n \to \infty}X(\varphi_n)$ in $L_{\bC}^2(\Omega)$, where $(\varphi_n) \subset \cD_{\bC}(\bR^{d+1})$ is such that $\|\varphi_n-\varphi\|_{0} \to 0$.

Using (\ref{def-X(varphi)}) and (\ref{cov-M}) (Appendix A), we have:
\begin{eqnarray*}E|X(\varphi_n)-Y|^2&=&E\left|\int_{\bR^{d+1}} [\cF \varphi_n(\tau,\xi)-\cF \phi_{\xi}(\tau)](-i\tau)\prod_{j=1}^{d}(-i \xi_j)\tM(d\tau,d\xi)\right|^2 \\
&=& \int_{\bR^{d+1}} |\cF \varphi_n(\tau,\xi)-\cF \phi_{\xi}(\tau)|^2\Pi(d\tau,d\xi)=\|\varphi_n-\varphi\|_{0}^2 \to 0,
\end{eqnarray*}
i.e. $Y=\lim_{n \to \infty}X(\varphi_n)$ in $L_{\bC}^2(\Omega)$. Hence, $Y=X(\varphi)$ a.s.
$\Box$

\vspace{3mm}

The next result gives a criterion for checking if the indicator function of a bounded Borel set $A$ (in particular, a bounded rectangle) is in $\cH$, and gives a representation for the stochastic integral of $1_{A}$.

\begin{theorem}
\label{indicators-in-H}
(i) If $A \in \cB_{b}(\bR^{d+1})$ is such that
\begin{equation}
\label{norm-1A-finite}\int_{\bR^{d+1}}|\cF 1_{A}(\tau,\xi)|^2 \Pi(d\tau,d\xi)<\infty,
\end{equation}
then $1_{A} \in \cH$ and the stochastic integral of $1_{A}$ with respect to $X$ is given by:
$$X(A):=X(1_{A})=\int_{\bR^{d+1}} \cF 1_{A}(\tau,\xi) (-i \tau)\prod_{j=1}^{d}(-i \xi_j)\tM(d\tau,d\xi)$$
(ii) Suppose that the measure $\widetilde{\Pi}$ satisfies the following condition: 
\begin{equation}
\label{cond-pi-tilde}\int_{\bR^{d+1}}\frac{\tau^2+|\xi|^2}{1+\tau^2+|\xi|^2} \widetilde{\Pi}(d\tau,d\xi)<\infty.
\end{equation}
Then (\ref{norm-1A-finite}) holds for any finite union $A$ of bounded rectangles in $\bR^{d+1}$. In particular, the stochastic integral of $1_{(0,t] \times (0,x]}$ with respect to $X$ is given by:
\begin{equation}
\label{def-X(t,x)}X(t,x):=X(1_{(0,t] \times (0,x]})=\int_{\bR^{d+1}} (e^{-i\tau t}-1)\prod_{j=1}^{d}(e^{-i \xi_j x_j}-1)\tM(d\tau,d\xi).
\end{equation}
\end{theorem}

\noindent {\bf Proof:} (i) Let $\varphi_n=1_{A} * \phi_n$ be the mollification of $1_{A}$, where $\phi_n(t,x)=n^{d+1}\phi(nt,nx)$, and $\phi \in \cD(\bR^{d+1})$ is such that
$\|\phi\|_{L^1(\bR^{d+1})}=1$. Then $\varphi_n \in \cD(\bR^{d+1})$, $\cF \varphi_n -\cF 1_{A}=\cF 1_{A}(\cF \phi_n-1)\to 0$ and $|\cF \phi_n-1| \leq 2$. By the dominated convergence theorem,
$$\int_{\bR^{d+1}}|\cF \varphi_n(\tau,\xi)-\cF 1_{A}(\tau,\xi)|^2 \Pi(d\tau,d\xi)\to 0.$$
This proves that $1_{A} \in \cH$. The second statement follows by an approximation argument, as in the last part of the proof of Theorem \ref{varphi-in-H}.

(ii) By linearity, it is enough to consider the case
when $A=(0,t] \times (0,x]$ with $t \in \bR, x \in \bR^{d}$. Note that (\ref{cond-pi-tilde}) is equivalent to $\int_{\bR^{d+1}} [1 \wedge (\tau^2+|\xi|^2)]\widetilde{\Pi}(d\tau,d\xi)<\infty$.
On the region $R=\{(\tau,\xi) \in \bR^{d+1}; \tau^2+|\xi|^2 \leq 1\}$, we use the fact that $|\cF 1_{A}| \leq |A|$, where $|A|$ is the Lebesgue measure of $A$, and hence $$\int_{R}|\cF 1_{A} (\tau,\xi)|^2 \Pi(d\tau,\xi) \leq |A|^2 \int_{R}\Pi(d\tau,\xi)<\infty.$$
 On the region $R^c$, we use the fact that $$|\cF 1_{A}(\tau,\xi)|=\left|\frac{e^{-i \tau t}-1}{-i \tau}\prod_{j=1}^{d}\frac{e^{-i \xi_j x_j}-1}{-i \xi_j}\right| \leq 2^{d+1}\left|\frac{1}{\tau \xi_1 \ldots \xi_d}\right|,$$
 and hence,
 $$\int_{R^c}|\cF 1_{A}(\tau,\xi)|^2 \Pi(d\tau,d\xi) \leq 2^{d+1}\int_{R^c}\widetilde{\Pi}(d\tau,d\xi)<\infty.$$
$\Box$

When condition (\ref{cond-pi-tilde}) holds, there is a one-to-one correspondence between the class of random fields $\{X(t,x)\}_{(t,x) \in \bR^{d+1}}$ defined by (\ref{def-X(t,x)}) and the class of random fields with stationary increments and spectral representation:
\begin{equation}
\label{def-X*(t,x)}X^*(t,x)=\int_{\bR^{d+1}} (e^{-i\tau t-i \xi \cdot x}-1)\tM(d\tau,d\xi).
\end{equation}
This correspondence is illustrated by the next two examples.

\begin{example}
\label{fBs-example}
{\rm The {\em fractional Brownian sheet} (fBs) with indices $H,H_1,\ldots,H_d \in (0,1)$ is a zero-mean Gaussian random field $\{X(t,x)\}_{(t,x) \in \bR^{d+1}}$ with covariance
$$E[X(t,x)X(s,y)]=R_{H}(t,s)\prod_{j=1}^{d}R_{H_j}(x_j,y_j),$$
where $R_{H}(t,s)=(|t|^{2H}+|s|^{2H}-|t-s|^{2H})/2$ is the covariance of the fBm of index $H \in (0,1)$.
 Since
$$R_{H}(t,s)=c_{H}\int_{\bR} (e^{-i \tau t}-1)(e^{i\tau s}-1)|\tau|^{-(2H+1)}d\tau,$$
with $c_{H}=\Gamma(2H+1)\sin(\pi H)/(2\pi)$, the fBs has the representation (\ref{def-X(t,x)}), where $\tM$ is Gaussian and $\widetilde{\Pi}(d\tau,d\xi)=c_{H}|\tau|^{-(2H+1)}\prod_{j=1}^{d}(c_{H_j}|\xi_j|^{-(2H_j+1)})$.
The measure $\widetilde{\Pi}$ satisfies (\ref{cond-pi-tilde}) if and only if $H+\sum_{j=1}^{d}H_j<1$. In this case, there is an associated random field with stationary increments defined by (\ref{def-X*(t,x)}), which is also called a ``fractional Brownian sheet'' by some authors (see \cite{clausel-vesel11}).
}
\end{example}

\begin{example}
\label{fBf-example}
{\rm A {\em fractional Brownian field} (fBf) on $\bR^{d+1}$ with index $H \in (0,1)$ is
a zero-mean Gaussian random field $\{X^*(t,x)\}_{(t,x) \in \bR^{d+1}}$ with covariance
$$E[X^*(t,x)X^*(s,y)]=\frac{1}{2}[(t^2+|x|^2)^{H}+(s^2+|y|^2)^{H}-(|t-s|^2+|x-y|^2)^{H}].$$
This random field has stationary increments and spectral representation (\ref{def-X*(t,x)}), where $\tM$ is Gaussian and $\widetilde{\Pi}(d\tau,d\xi)=(\tau^2+|\xi|^2)^{-(2H+d+1)/2}$. The measure $\widetilde{\Pi}$ satisfies condition (\ref{cond-pi-tilde}) for any $H \in (0,1)$.
To this process, one may associate a random field $\{X(t,x)\}$ defined by (\ref{def-X(t,x)}).


}
\end{example}

An interesting case arises when
\begin{equation}
\label{Pi-product}
\Pi=\nu \times \mu,
\end{equation}
 where $\nu$ is a symmetric measure on $\bR$ satisfying (\ref{nu-cond}), and $\mu$ is a symmetric tempered measure on $\bR^d$.
Using the same mollification argument as in the proof of Theorem \ref{indicators-in-H}.(i), it follows that $1_{[0,t]} \varphi \in \cH$ for any $t \in \bR,\varphi \in \cD_{\bC}(\bR^d)$. To see this, note that
$$\int_{\bR^{d+1}} |\cF (1_{[0,t]}\varphi)(\tau,\xi)|^2\Pi(d\tau,d\xi)=\int_{\bR}\frac{|e^{-i \tau t}-1|^2}{ \tau^2}\nu(d\tau)\int_{\bR^d}|\cF \varphi(\xi)|^2\mu(d\xi)<\infty.$$

For any $t,s \in \bR$ and $\varphi,\psi \in \cD_{\bC}(\bR^d)$, we have:
$$\langle 1_{[0,t]}\varphi,1_{[0,s]}\psi \rangle_{\cH}=R(t,s)J(\varphi,\psi),$$
where $R$ is given by (\ref{time-covariance}), being the covariance of a process $Z=(Z_t)_{t \in \bR}$ with stationary increments and spectral measure $\tilde{\nu}(d\tau)=\tau^{-2}\nu(d\tau)$, and
$$J(\varphi,\psi)=\int_{\bR^d} \cF \varphi(\xi)\overline{\cF \psi(\xi)}\mu(d\xi).$$

In this case, $\cH$ is the completion of $\cE$ with respect to the inner
product $\langle \cdot, \cdot \rangle_{\cH}$, where $\cE$ is the set of complex linear combinations of functions of the form $1_{[0,t]} \varphi$, with $t \in \bR$ and $\varphi \in \cD_{\bC}(\bR^d)$.
To prove this, we use that fact that a function in $\cD_{\bC}(\bR^{d+1})$ can be approximated by a function of the form $\varphi(t,x)=\varphi_1(t)\varphi_2(x)$ with
$\varphi_1\in \cD_{\bC}(\bR),\varphi_2 \in \cD_{\bC}(\bR^d)$, and a function in $\cD_{\bC}(\bR)$ can approximated in the norm $\|\cdot\|_{\Lambda}$ by a complex elementary function on $\bR$ (see Remark \ref{Lambda-rem}).

The stochastic integral of $1_{[0,t]}\varphi$ with respect to $X$ is given by:
$$X_t(\varphi):=X(1_{[0,t]}\varphi)=\int_{\bR^{d+1}} (e^{-i\tau t}-1)\cF \varphi(\xi)  \prod_{j=1}^{d}(-i \xi_j)\tM(d\tau,d\xi).$$

\section{Equations with harmonizable noise}
\label{equations-section}

As in the previous section, let $\{X(\varphi);\varphi \in \cD_{\bC}(\bR^{d+1})\}$ be a harmonizable noise with driving measure $\tM$ and covariance measure $\Pi$.
We consider the equation:
\begin{equation}
\label{spde} L u(t,x)=\dot{X}(t,x) \quad t>0, x \in \bR^d
\end{equation}
with some deterministic initial conditions, where $L$ is a second order partial differential operator. Let $w$ be the
solution of the equation $Lu=0$ with the same initial conditions
as (\ref{spde}), and $G$ be
the fundamental solution of $Lu=0$.

The following definition introduces the concept of solution.

\begin{definition}
{\rm The process $\{u(t,x);t \geq 0, x\in \bR^d\}$ defined by
\begin{equation}
\label{def-sol}u(t,x)=w(t,x)+\int_0^t
\int_{\bR^d}G(t-s,x-y)X(ds,dy)
\end{equation}
is called a {\em random field solution} of (\ref{spde}), provided
that the stochastic integral on the right-hand side of
(\ref{def-sol}) is well-defined, i.e. for any $t>0$ and $x \in
\bR^d$
$$g_{tx}:=1_{[0,t]}G(t-\cdot,x-\cdot) \in \cH.$$}
\end{definition}

The following result gives some general conditions which ensure the existence of a random field solution, as a direct consequence of Theorem
\ref{varphi-in-H}.

\begin{theorem}
\label{existence-sol-th}
Assume that $G(t,\cdot) \in \cS'(\bR^d)$, such that $\cF
G(t,\cdot)$ is a function for any $t>0$, $(t,\xi) \mapsto \cF G(t,\cdot)(\xi)=:H_{\xi}(t)$
is measurable on $\bR_{+} \times \bR^d$, and
$$\int_{0}^{t}|\cF G(s,\cdot)(\xi)|ds<\infty, \quad \mbox{for any} \ t>0,\xi \in \bR^d.$$

\noindent (i) Equation (\ref{spde}) has a random field solution $\{u(t,x)\}_{(t,x) \in \bR^{d+1}}$ if and only if
\begin{equation}
\label{It-is-finite} I_t:=\int_{\bR^{d+1}} |\cF_{0,t}
H_{\xi}(\tau)|^2 \Pi(d\tau,d\xi)<\infty \quad \mbox{for all} \ t>0,
\end{equation}
where
$\cF_{0,t} H_{\xi} (\tau)=\int_{0}^{t}e^{-i \tau s}H_{\xi}(s)ds$,
$\tau \in \bR$. In this case, $E(u(t,x))=w(t,x)$, ${\rm Var}(u(t,x))=I_t$
and the solution is given by:
$$u(t,x)=w(t,x)+\int_{\bR^{d+1}}e^{-i\tau t-i \xi \cdot x}\overline{\cF_{0,t}H_{\xi}(\tau)}(-i \tau) \prod_{j=1}^{d}(-i \xi_j) \tM(d\tau,d\xi).$$

\noindent (ii) Assume zero initial conditions. Then, for any $t,s>0$ and $x,y \in \bR^d$,
$$E[u(t,x)u(s,y)]=\int_{\bR^{d+1}}e^{-i \tau (t-s)}e^{-i \xi \cdot (x-y)} \overline{\cF_{0,t}H_{\xi}(\tau)} \cF_{0,s} H_{\xi}(\tau) \Pi(d\tau,d\xi).$$
In particular, for any $t>0$ fixed,
$$E[u(t,x)u(t,y)]=\int_{\bR^{d}}e^{-i \xi \cdot (x-y)}\mu_t(d\xi).$$
where $\mu_t(\cdot)=\Pi_t(\bR \times \cdot)$ and $\Pi_t(d\tau,d\xi)=|\cF_{0,t}H_{\xi}(\tau)|^2 \Pi(d\tau,d\xi)$, i.e. $\{u(t,x)\}_{x \in \bR^d}$ is a zero-mean stationary random field on $\bR^d$ with spectral measure $\mu_t$.
\end{theorem}

\noindent {\bf Proof:} (i) We apply Theorem \ref{varphi-in-H} to
$\varphi=g_{tx}$. Note that the function
$$\phi_{\xi}(s):=\cF g_{tx}(s,\cdot)(\xi)= 1_{[0,t]}(s)e^{-i \xi \cdot x}\overline{\cF G(t-s,\cdot)(\xi)}
=1_{[0,t]}(s)e^{-i \xi \cdot x} \overline{H_{\xi}(t-s)}$$
satisfies conditions (i) and (ii) of Theorem \ref{varphi-in-H}.

Let $\cF \phi_{\xi}$ be the Fourier transform of the function $s
\mapsto \phi_{\xi}(s)$. Then,
\begin{eqnarray}
\nonumber
\cF \phi_{\xi}(\tau)&=&\int_{\bR} e^{-i \tau s} \phi_{\xi}(s)ds = e^{-i \xi \cdot x} \int_0^t e^{-i \tau s} \overline{H_{\xi}(t-s)}ds \\
\label{Fourier-phi-xi}
&=& e^{-i \xi \cdot x} e^{-i \tau t} \int_0^t e^{i \tau s}
\overline{H_{\xi}(s)}ds= e^{-i \xi \cdot x} e^{-i \tau t}
\overline{\cF_{0,t} H_{\xi}(\tau)}.
\end{eqnarray}

By definition $u(t,x)=w(t,x)+X(g_{tx})$. Hence, $E(u(t,x))=w(t,x)$ and
${\rm Var}(u(t,x))=E|X(g_{tx})|^2=\|g_{tx}\|_{\cH}^2=I_t$.

(ii) Using the isometry property of the stochastic integral $\varphi \mapsto X(\varphi)$ and the fact that $\|\varphi\|_{\cH}=\|\varphi\|_{0}$, where $\|\cdot \|_{0}$ is defined by relation (\ref{norm-varphi-0-finite}), we obtain:
\begin{eqnarray*}
E[u(t,x)u(s,y)]&=& E[X(g_{tx})X(g_{sy})]=\langle g_{tx}, g_{sy}\rangle_{\cH}
= \langle g_{tx}, g_{sy}\rangle_{0} \\
&=& \int_{\bR^{d+1}}\cF \phi_{\xi}(\tau) \overline{\cF \psi_{\xi}(\tau)}\Pi(d\tau,d\xi),
\end{eqnarray*}
where $\phi_{\xi}(u)=\cF g_{tx}(u,\cdot)(\xi)$ and $\psi_{\xi}(u)=\cF g_{sy}(u,\cdot)(\xi)$. The conclusion follows using (\ref{Fourier-phi-xi}). $\Box$

\section{A parabolic equation}
\label{parabolic-section}

In this section, we assume that (\ref{Pi-product}) holds.
We consider the equation:
\begin{equation}
\label{parabolic-eq} \frac{\partial u}{\partial t}(t,x)+
(-\Delta)^{\beta/2}u(t,x)=\dot{X}(t,x), \quad t>0,x \in \bR^d,
\end{equation}
for some $\beta>0$. In this case, $$H_{\xi}(t)=\cF G
(t,\cdot)=\exp \{-t \Psi(\xi)\},$$ where
$\Psi(\xi)=c_{\beta}|\xi|^{\beta}$ and $c_{\beta}>0$ is a constant.
By Theorem \ref{existence-sol-th}, it suffices to
check that (\ref{It-is-finite}) holds. For this, we write
$I_t=\int_{\bR^d} N_{t}(\xi)\mu(d\xi)$, where
\begin{equation}
\label{def-Nt} N_t(\xi)=\int_{\bR} |\cF_{0,t} H_{\xi}(\tau)|^2
\nu(d\tau).
\end{equation}

We proceed to the calculation of $\cF_{0,t} H_{\xi}(\tau)$. For this,
note that if $\varphi(x)=e^{-ax}$ with $a \in \bR$, then
$$\cF_{0,t} \varphi(\tau)=\int_0^t e^{-(a+i \tau) x}dx=\frac{1-e^{-(a+i\tau)t}}{a+i\tau}=\frac{1-e^{-at} \cos(\tau t)+i e^{-at} \sin (\tau t)}{a+i \tau}$$
and
\begin{equation}
\label{Fourier-expo} |\cF_{0,t} \varphi
(\tau)|^2=\frac{1}{\tau^2+a^2} \{\sin^2(\tau t)+[e^{-at}-\cos(\tau
t)]^2 \}.
\end{equation}

Applying (\ref{Fourier-expo}) with $a=\Psi(\xi)$, we obtain:
\begin{equation}
\label{Fourier-Nt-par} |\cF_{0,t}
H_{\xi}(\tau)|^2=\frac{1}{\tau^2+\Psi(\xi)^2}\{\sin^2(\tau
t)+[e^{-t\Psi(\xi)}- \cos(\tau t)]^2\}.
\end{equation}

Using (\ref{def-Nt}) and (\ref{Fourier-Nt-par}), we obtain:
$$N_t(\xi)=\int_{\bR^d}\frac{1}{\tau^2+\Psi(\xi)^2}\{\sin^2(\tau t)+[e^{-t\Psi(\xi)}- \cos(\tau t)]^2\}\nu(d\tau).$$

To identify the necessary and sufficient condition for the existence of a solution, it suffices to obtain some suitable estimates for $N_t(\xi)$. This goal is achieved in the next theorem.


\begin{theorem}
\label{parabolic-estimates} For any $t>0$ and $\xi \in \bR^d$, we
have:
$$C_t^{(2)} \int_{\bR^d} \frac{1}{\tau^2+[1+\Psi(\xi)]^2} \nu(d\tau) \leq N_t(\xi) \leq C_t^{(1)}
\int_{\bR^d} \frac{1}{\tau^2+[1+\Psi(\xi)]^2} \nu(d\tau),$$ where
$C_t^{(1)}$ and $C_t^{(2)}$ are some positive constants depending
on $t$. Consequently, equation (\ref{parabolic-eq}) has a random
field solution if and only if (\ref{parabolic-cond}) holds.
\end{theorem}

\begin{remark}
{\rm In the case of Example \ref{Riesz-ex} and Example
\ref{Bessel-ex} with $\gamma \in (-1,1)$, condition
(\ref{parabolic-cond}) becomes: (use Lemma \ref{lemmaA}, Appendix
B with $a=1+\Psi(\xi)$)
\begin{equation}
\label{cond-par-Riesz}
\int_{\bR^d} \left( \frac{1}{1+\Psi(\xi)}\right)^{1+\gamma}\mu(d\xi)<\infty.
\end{equation}
When applied to Example \ref{Riesz-ex} with the parametrization
$\gamma=2H-1,H \in (0,1)$, Theorem \ref{parabolic-estimates}
becomes an extension of Theorem 3.4 of \cite{B11} to the case
$H<1/2$. }
\end{remark}

\noindent {\bf Proof of Theorem \ref{parabolic-estimates}:} {\em
(a) We first show the upper bound.} We claim that:
\begin{equation}
\label{Nt-heat-bound1} N_t(\xi) \leq 2K \max(t^2,5), \quad
\mbox{for all} \ \xi \in \bR^d.
\end{equation}

To prove (\ref{Nt-heat-bound1}), we split the integral on the
right-hand side of (\ref{def-Nt}) into the regions $|\tau| \leq 1$
and $|\tau| \geq 1$. We denote the two integrals by
$N_t^{(1)}(\xi)$ and $N_t^{(2)}(\xi)$. For $N_t^{(1)}(\xi)$, we
use the fact that: $$|\cF_{0,t} H_{\xi}(\tau)|=\left|\int_0^t
e^{-i \tau s} e^{-s \Psi(\xi)}ds\right| \leq \int_0^t e^{-s
\Psi(\xi)}ds \leq t,$$ and hence
\begin{equation}
\label{I1-bound} N_t^{(1)}(\xi)=\int_{|\tau| \leq 1}|\cF_{0,t}
H_{\xi}(\tau)|^2d\nu(\tau) \leq t^2 \int_{|\tau| \leq 1}
\nu(d\tau) \leq t^2 \int_{|\tau| \leq 1}\frac{2}{1+\tau^2}
\nu(d\tau).
\end{equation}
For $N_t^{(2)}(\xi)$, we use the fact that:
\begin{equation}
\label{5-sin-cos-bound} \sin^2(\tau t)+[e^{-t \Psi(\xi)}-\cos(\tau
t)]^2 \leq 5,
\end{equation}
and hence
\begin{equation}
\label{I2-bound}N_t^{(2)}(\xi) \leq 5 \int_{|\tau| \geq 1}
\frac{1}{\tau^2+\Psi(\xi)^2} \leq 5 \int_{|\tau| \geq 1}
\frac{1}{\tau^2}\nu(d\tau) \leq 5 \int_{|\tau| \geq
1}\frac{2}{\tau^2+1}\nu(d\tau).
\end{equation}
Relation (\ref{Nt-heat-bound1}) follows from (\ref{I1-bound}) and
(\ref{I2-bound}).
On the other hand, by (\ref{5-sin-cos-bound}),
\begin{equation}
\label{Nt-heat-bound4} N_t(\xi) \leq 5 \int_{\bR}
\frac{1}{\tau^2+\Psi(\xi)^2} \nu(d\tau), \quad \mbox{for all} \
\xi \in \bR^d.
\end{equation}

{\em (i) Assume $\Psi(\xi) \leq 1$.} Then $\tau^2+[1+\Psi(\xi)]^2
\leq \tau^2+4 \leq 4(\tau^2+1)$, and hence
\begin{equation}
\label{Nt-heat-bound2}\int_{\bR}
\frac{1}{\tau^2+[1+\Psi(\xi)]^2}\nu(d\tau) \geq \frac{1}{4}
\int_{\bR}\frac{1}{1+\tau^2}\nu(d\tau)=\frac{1}{4}K, \quad
\mbox{if} \ \Psi(\xi) \leq 1.
\end{equation}

From (\ref{Nt-heat-bound1}) and (\ref{Nt-heat-bound2}), we infer
that:
\begin{equation}
\label{Nt-heat-bound3}N_t(\xi) \leq 8 \max(t^2,5)
\int_{\bR}\frac{1}{\tau^2+[1+\Psi(\xi)]^2}\nu(d\tau) \quad
\mbox{if} \ \Psi(\xi) \leq 1.
\end{equation}

{\em (ii) Assume $\Psi(\xi) \geq 1$.} Then $4[\tau^2+\Psi(\xi)]^2
\geq 4\tau^2 +[1+\Psi(\xi)]^2 \geq \tau^2+[1+\Psi(\xi)]^2$, and
hence
\begin{equation}
\label{Nt-heat-bound5} \int_{\bR}
\frac{1}{\tau^2+\Psi(\xi)^2}\nu(d\tau) \leq 4 \int_{\bR}
\frac{1}{\tau^2+[1+\Psi(\xi)]^2}\nu(d\tau), \quad \mbox{if} \
\Psi(\xi) \geq 1.
\end{equation}

From (\ref{Nt-heat-bound4}) and (\ref{Nt-heat-bound5}), we infer
that:
\begin{equation}
\label{Nt-heat-bound6} N_t(\xi) \leq 20 \int_{\bR}
\frac{1}{\tau^2+[1+\Psi(\xi)]^2}\nu(d\tau), \quad \mbox{if} \
\Psi(\xi) \geq 1.
\end{equation}

The upper bound follows from (\ref{Nt-heat-bound3}) and
(\ref{Nt-heat-bound6}), with
$C_t^{(1)}=8 \max(t^2,5)$.

{\em (b) We now show the lower bound.} Let $\rho_t>0$ be a
constant depending on $t$, which will be specified below.

{\em (i) Assume $\Psi(\xi) \leq \rho_t$.}  We use the fact that
$\sin x$ is decreasing on $[\frac{\pi}{2},\pi]$. Let
$\frac{\pi}{2}<c<d<\pi$ be fixed. We have:
\begin{eqnarray}
\nonumber N_t(\xi) & \geq & \int_{\bR}
\frac{1}{\tau^2+\Psi(\xi)^2} \sin^2(\tau t) \nu(d\tau)
\geq  \int_{c \leq |\tau t| \leq d} \frac{1}{\tau^2+\Psi(\xi)^2} \sin^2(\tau t) \nu(d\tau) \\
\nonumber
& \geq & \sin^2 d \int_{c \leq |\tau t| \leq d}\frac{1}{\tau^2+\Psi(\xi)^2}\nu(d\tau) \\
\label{Nt-heat-LB1} & \geq & \frac{\sin^2
d}{d^2/t^2+\rho_t^2}\int_{c/t \leq |\tau| \leq
d/t}\nu(d\tau)=:A(t).
\end{eqnarray}

On the other hand,
\begin{equation}
\label{Nt-heat-LB2} K=\int_{\bR} \frac{1}{\tau^2+1}\nu(d\tau) \geq
\int_{\bR}\frac{1}{\tau^2+[1+\Psi(\xi)]^2} \nu(d\tau), \quad
\mbox{for all} \ \xi \in \bR^d.
\end{equation}

From (\ref{Nt-heat-LB1}) and (\ref{Nt-heat-LB2}), we obtain:
\begin{equation}
\label{Nt-heat-LB3} N_t(\xi) \geq \frac{A(t)}{K}
\int_{\bR}\frac{1}{\tau^2+[1+\Psi(\xi)]^2} \nu(d\tau).
\end{equation}

{\em (ii) Assume $\Psi(\xi) \geq \rho_t$.} For any $x \geq 0$, let
$$f_{\tau,t}(x)=\sin^2(\tau t)+[e^{-tx}-\cos(\tau t)]^2=1-[2e^{-tx} \cos(\tau t)-e^{-2tx}].$$
Note that $f_{t,\tau}(x) \to 1$ as $x \to \infty$, uniformly in
$\tau \in \bR$. Hence, there exists $\rho_t>0$ such that
\begin{equation}
\label{Nt-heat-LB4} f_{\tau,t}(x) \geq \frac{1}{2}, \quad
\mbox{for all $x \geq \rho_t$ and $\tau \in \bR$}.
\end{equation}
(More precisely, $|2e^{-tx} \cos(\tau t)-e^{-2tx}| \leq
2e^{-tx}+e^{-2tx} \leq 2c+c^2$, if $e^{-tx} \leq c$. Choose $c \in
(0,1)$ such that $2c+c^2 \leq 1/2$, e.g. $c=0.1$. Then
$f_{t,\tau}(x) \geq 1/2$ for any $x \geq (-\ln c)/t=:\rho_t$.)
Using (\ref{Nt-heat-LB4}), we infer that:
\begin{equation}
\label{Nt-heat-LB5} N_t(\xi) \geq \frac{1}{2}
\int_{\bR}\frac{1}{\tau^2+\Psi(\xi)^2}\nu(d\tau) \geq \frac{1}{2}
\int_{\bR}\frac{1}{\tau^2+[1+\Psi(\xi)]^2}\nu(d\tau), \quad
\mbox{if} \ \Psi(\xi) \geq \rho_t.
\end{equation}
The lower bound follows from (\ref{Nt-heat-LB3}) and
(\ref{Nt-heat-LB5}), with
$C_{t}^{(2)}=\min \left\{\frac{A(t)}{K},\frac{1}{2} \right\}$.
$\Box$

\begin{remark}
{\rm The results of this section are also valid for the equation $$\frac{\partial u}{\partial t}(t,x)-\cL
u(t,x)=\dot{X}(t,x), \quad t>0,x \in \bR^d,$$
where $\cL$ is the $L^2(\bR^d)$-generator of a $d$-dimensional
L\'evy process $(X_t)_{t \geq 0}$ with (real-valued) characteristic exponent
$\Psi(\xi)$ (as in \cite{FK10}). Assuming that $X_t$ has density
$p_{t}$, we see that $G(t,x)=p_t(-x)$ and
$$H_{\xi}(t)=\cF G(t,\cdot)(\xi)=\int_{\bR^d}e^{i \xi \cdot x}p_t(x)dx=E(e^{i \xi \cdot X_t})=\exp\{-t \Psi(\xi)\}.$$
}
\end{remark}

\section{A hyperbolic equation}
\label{hyperbolic-section}

As in Section \ref{parabolic-section}, we assume that (\ref{Pi-product}) holds. We consider the equation:
\begin{equation}
\label{hyperbolic-eq} \frac{\partial^2 u}{\partial t^2}(t,x)+(-\Delta)^{\beta/2}
u(t,x)=\dot{X}(t,x), \quad t>0, x \in \bR^d,
\end{equation}
for some $\beta>0$. In this case,
$$H_{\xi}(t)=\cF G(t,\cdot)(\xi)=\frac{\sin(t \sqrt{\Psi(\xi)})}{\sqrt{\Psi(\xi)}}$$
where $\Psi(\xi)=|\xi|^{\beta}$ (see (16) of
\cite{dalang-sanzsole05}).

Following the same steps as in the parabolic case,
it suffices to find appropriate bounds for $I_t=\int_{\bR^d}N_t(\xi)\mu(d\xi)$, where $N_t(\xi)$ is given by (\ref{def-Nt}).
For the calculation of $\cF_{0,t} H_{\xi}(\tau)$, we use the change of variable
$r=s\sqrt{\Psi(\xi)}$:
\begin{eqnarray}
\label{Fourier-Nt-hyp1}
\cF_{0,t} H_{\xi}(\tau)&=&\frac{1}{\sqrt{\Psi(\xi)}}\int_{0}^{t} e^{-i\tau s} \sin(s \sqrt{\Psi(\xi)})ds \\
\nonumber
&=& \frac{1}{\Psi(\xi)} \int_{0}^{t\sqrt{\Psi(\xi)}} e^{-i\tau r/\sqrt{\Psi(\xi)}} \sin r dr \\
\nonumber &=& \frac{1}{\Psi(\xi)} \cF_{0,T} \sin
\left(\frac{\tau}{\sqrt{\Psi(\xi)}} \right), \quad \mbox{where} \
T=t \sqrt{\Psi(\xi)}.
\end{eqnarray}

An elementary calculation shows that: (see the proof of Lemma B.1
of \cite{BT10})
\begin{equation}
\label{Fourier-sin} |\cF_{0,T} \sin
(\tau)|^2=\frac{1}{(\tau^2-1)^2}[f_{T}^2(\tau)+g_T^2(\tau)].
\end{equation}
where
$f_{T}(\tau)=\sin(\tau T)-\tau \sin T$ and $g_{T}(\tau)=\cos(\tau T)-\cos T$.
Hence,
\begin{equation}
\label{Fourier-Nt-hyp2} |\cF_{0,t} H_{\xi}(\tau)|^2=
\frac{[\sin(\tau t)-\frac{\tau}{\sqrt{\Psi(\xi)}}
\sin(t\sqrt{\Psi(\xi)})]^2+[\cos(\tau t)-\cos(t
\sqrt{\Psi(\xi)})]^2}{(\tau^2-\Psi(\xi))^2}.
\end{equation}

Using (\ref{def-Nt}) and (\ref{Fourier-Nt-hyp2}), we obtain:
$$N_t(\xi)=\int_{\bR} \frac{[\sin(\tau t)-\frac{\tau}{\sqrt{\Psi(\xi)}} \sin(t\sqrt{\Psi(\xi)})]^2+[\cos(\tau t)-\cos(t \sqrt{\Psi(\xi)})]^2}{(\tau^2-\Psi(\xi))^2} \nu(d\tau).$$

We assume that $$\nu(d \tau)=\eta(|\tau|)d\tau,$$ where the
function $\eta$ satisfies the following condition:
\begin{eqnarray*}
(C) & & \mbox{for any $\lambda>0$ there exists $C_{\lambda}>0$
such that $\eta(\lambda \tau) \leq C_{\lambda} \eta(\tau) \
\forall \tau>0$.}
\end{eqnarray*}

In addition, we assume that $\eta$ satisfies either (C1) or (C2)
below:
\begin{eqnarray*}
(C1) & & \mbox{$\eta$ is non-increasing on $(0,\infty)$, and for any $K>0$ there exists $D_K>0$} \\ & & \mbox{such that $\int_0^a \eta(\tau)d\tau \leq D_{K} a \eta(a)$ for any $a \geq K$ }\\
(C2) & & \mbox{$\eta$ is non-decreasing on $(0,\infty)$, and for
any $K>0$ there exists $D_K>0$} \\ & & \mbox{such that
$\int_a^{\infty} \tau^{-2}\eta(\tau)d\tau \leq D_{K} a^{-1}
\eta(a)$ for any $a \geq K$. }
\end{eqnarray*}

\begin{remark}
{\rm In Examples \ref{Riesz-ex} and \ref{Bessel-ex}, condition (C)
holds for any $\gamma \in (-1,1)$, respectively for any
$\gamma>-1$. In these two examples, condition (C1) holds if
$\gamma \in [0,1)$, whereas condition (C2) holds if $\gamma \in
(-1,0)$ (see Appendix C). } 
\end{remark}

Under these conditions, we obtain the following result.

\begin{theorem}
\label{main-th-hyp} Suppose that $\eta$ satisfies (C). In
addition, suppose that $\eta$ satisfies either (C1) or (C2). Then
for any $t>0$ and $\xi \in \bR^d$,
$$ D_t^{(2)}N(\xi) \leq N_t(\xi) \leq D_t^{(1)} N(\xi),$$
where $D_t^{(1)}$ and $D_t^{(2)}$ are positive constants depending
on $t$, and
$$N(\xi)=\frac{1}{\sqrt{1+\Psi(\xi)}} \int_{\bR} \frac{1}{(|\tau|+\sqrt{1+\Psi(\xi)})^2}\nu(d\tau).$$
Consequently, equation (\ref{hyperbolic-eq}) has a random field
solution if and only if
\begin{equation}
\label{hyperbolic-cond} \int_{\bR^d}\frac{1}{\sqrt{1+\Psi(\xi)}}
\int_{\bR} \frac{1}{(|\tau|+\sqrt{1+\Psi(\xi)})^2}\nu(d\tau)
\mu(d\xi)<\infty.
\end{equation}
\end{theorem}

\begin{remark}
{\rm Note that condition (\ref{hyperbolic-cond}) is
equivalent to (\ref{hyperbolic-cond2}). This follows using the fact that $\tau^2+b^2 \leq (|\tau|+b)^2 \leq 2(\tau^2+b^2)$ with
$b=\sqrt{1+\Psi(\xi)}$.
}
\end{remark}

The proof of the theorem follows from the lemmas below. The first
two lemmas treat the upper bound.

\begin{lemma}
\label{hyp-UB-lemma1} Assume that $\eta$ satisfies (C) and is
either non-increasing on $(0,\infty)$, or non-decreasing on
$(0,\infty)$. Then for any $t>0$ and $\xi \in \bR^d$,
$$N_t(\xi) \leq C t \frac{1}{\sqrt{\Psi(\xi)}}\int_{\bR} \frac{1}{(|\tau|+\sqrt{\Psi(\xi)})^2}\nu(d\tau),$$
where $C$ is a positive constant.
\end{lemma}

\noindent {\bf Proof:} Using the notation $a=\sqrt{\Psi(\xi)}$, we have:
\begin{equation}
\label{def-Nt-hyp}
N_t(\xi)=\int_{\bR}\frac{1}{(\tau^2-a^2)^2}\left[f_{ta}^2\left(\frac{\tau}{a}\right)+
g_{ta}^2\left(
\frac{\tau}{a}\right)\right]\nu(d\tau).
\end{equation}

We denote by $N_t^{(1)}(\xi)$, $N_t^{(2)}(\xi)$ and
$N_t^{(3)}(\xi)$ the integrals over the regions $R_1=\{|\tau| \leq
a/2\}$, $R_2=\{|\tau| \geq 3a/2\}$, respectively $R_3=\{a/2 \leq
|\tau| \leq 3a/2\}$.

{\em (i) We first treat $N_t^{(1)}(\xi)$ and $N_t^{(2)}(\xi)$.}
For this, we use the inequality:
\begin{equation}
\label{basic-ineq} f_T^2(\tau)+g_T^2(\tau) \leq 2T (1+|\tau|)^2
\quad \mbox{for any} \quad \tau \in \bR.
\end{equation}
(To see this, note that $|f_T(\tau)| \leq 1+|\tau|$ and
$|f_T(\tau)| \leq 2 T|\tau|$, since $|\sin x| \leq |x|$. Hence,
$f_T^2(\tau) \leq 2 T |\tau|(1+|\tau|)$. Similarly, $|g_T(\tau)|
\leq 2$ and $|g_T(\tau)| \leq T(1+|\tau|)$, since $|1-\cos x| \leq
|x|$. Hence, $g_T^2(\tau) \leq 2T (1+|\tau|)$.)

Using (\ref{basic-ineq}), we obtain:
\begin{equation}
\label{important-ineq} f_{ta}^{2} \left(\frac{\tau}{a}
\right)+g_{ta}^2\left( \frac{\tau}{a}\right) \leq  2ta \left(
1+\frac{|\tau|}{a}\right)^2=\frac{2t}{a}(|\tau|+a)^2.
\end{equation}
Therefore,
$$N_t^{(i)}(\xi) \leq  \frac{2t}{a} \int_{R_i} \frac{1}{(|\tau|-a)^2}\nu(d\tau), \quad \mbox{for} \quad i=1,2.$$
When $\tau \in R_1$, both $(|\tau|-a)^{-2}$ and $(|\tau|+a)^{-2}$
behave as $a^{-2}$, since
$$\frac{1}{a^2} \leq \frac{1}{(|\tau|-a)^2} \leq \frac{4}{a^2} \quad \mbox{and} \quad
\frac{4}{9a^2} \leq \frac{1}{(|\tau|+a)^2} \leq \frac{1}{a^2}.$$
When $\tau \in R_2$, both $(|\tau|-a)^{-2}$ and $(|\tau|+a)^{-2}$
behave as $\tau^{-2}$, since
$$\frac{1}{\tau^2} \leq \frac{1}{(|\tau|-a)^2} \leq \frac{9}{\tau^2} \quad \mbox{and} \quad
\frac{9}{25 \tau^2} \leq \frac{1}{(|\tau|+a)^2} \leq
\frac{1}{\tau^2}.$$

Letting $C^{(1)}=18$ and $C^{(2)}=50$, we obtain
$$N_t^{(i)}(\xi) \leq  C^{(i)} \frac{t}{a} \int_{R_i} \frac{1}{(|\tau|+a)^2}\nu(d\tau), \quad \mbox{for} \quad i=1,2.$$

{\em (ii) We now treat $N_t^{(3)}(\xi)$.} From
(\ref{Fourier-Nt-hyp1}) and (\ref{Fourier-Nt-hyp2}), we see that:
\begin{equation}
\label{Fourier-sin2} \frac{1}{(\tau^2-a^2)^2} \left[f_{ta}^2
\left(\frac{\tau}{a}\right)+g_{ta}^2\left(\frac{\tau}{a}\right)\right]
=\frac{1}{a^2}|\cF_{0,t} \sin (a \ \cdot)(\tau)|^2.
\end{equation}

Assume first that $\eta$ is non-increasing on $(0,\infty)$. Then
\begin{eqnarray}
\nonumber
N_t^{(3)}(\xi)& =& \frac{1}{a^2} \int_{R_3} |\cF_{0,t} \sin (a \ \cdot)(\tau)|^2 \eta(|\tau|)d\tau \\
\nonumber
& \leq & \frac{\eta(a/2)}{a^2}\int_{\bR} |\cF_{0,t} \sin (a \ \cdot)(\tau)|^2 d\tau \\
\label{UB-Nt3-step1} &=& \frac{\eta(a/2)}{a^2} 2\pi \int_{0}^{t}
\sin^2(as)ds \leq 2 \pi t \frac{\eta(a/2)}{a^2},
\end{eqnarray}
where we used Plancherel's theorem for the second equality above.
On the other hand, when $\tau \in R_3$, $(|\tau|+a)^{-2}$ behaves
as $a^{-2}$, since:
$$\frac{4}{25 a^2} \leq \frac{1}{(|\tau|+a)^2} \leq \frac{4}{9a^2}.$$
Using the fact that $\eta$ is non-increasing, we obtain:
\begin{equation}
\label{UB-Nt3-step2} \frac{1}{a}\int_{R_3}
\frac{\eta(|\tau|)}{(|\tau|+a)^2} d\tau \geq \frac{1}{a}\cdot
\frac{4}{25 a^2}\eta(3a/2) \int_{R_3}d\tau=\frac{8}{25}
\frac{\eta(3a/2)}{a^2}.
\end{equation}

From condition (C), there exists a constant $C_{1/3}>0$ such that:
\begin{equation}
\label{UB-Nt3-step3} \eta(a/2) \leq C_{1/3} \eta(3a/2).
\end{equation}

From (\ref{UB-Nt3-step1}), (\ref{UB-Nt3-step2}) and
(\ref{UB-Nt3-step3}), it follows that:
$$N_t^{(3)}(\xi) \leq 2 \pi t \frac{25}{8}C_{1/3}\frac{1}{a}\int_{R_3} \frac{1}{(|\tau|+a)^2}\nu(d\tau).$$

When $\eta$ is non-decreasing, the argument is similar. In this
case,
$$N_t^{(3)}(\xi) \leq 2\pi t \frac{\eta(3a/2)}{a^2} \leq 2\pi t C_3 \frac{\eta(a/2)}{a^2} \leq 2 \pi t C_3 \frac{25}{8} \frac{1}{a}\int_{R_3} \frac{1}{(|\tau|+a)^2}\nu(d\tau).$$
$\Box$

\begin{lemma}
\label{hyp-UB-lemma2} Assume that $\eta$ satisfies (C) and is
either non-increasing on $(0,\infty)$, or non-decreasing on
$(0,\infty)$. Then for any $t>0$ and $\xi \in \bR^d$,
$$N_t(\xi) \leq D_t^{(1)} \frac{1}{\sqrt{1+\Psi(\xi)}}\int_{\bR} \frac{1}{(|\tau|+\sqrt{1+\Psi(\xi)})^2}\nu(d\tau),$$
where $D_t^{(1)}$ is a positive constant depending on $t$.
\end{lemma}

\noindent {\bf Proof:} We denote $a=\sqrt{\Psi(\xi)}$,
$b=\sqrt{1+\Psi(\xi)}$ and
$$I(x)=\int_{\bR}\frac{1}{(|\tau|+x)^2}\nu(d\tau), \quad \mbox{for any} \ x>0.$$

{\em (i) Assume that $a \leq 1$.} We first show that:
\begin{equation}
\label{Nt-wave-UB1} N_t(\xi) \leq 3K
\max\{\frac{t^2}{4},8(t^2+\frac{1}{2})\}.
\end{equation}

Using (\ref{def-Nt-hyp}), we denote by $N_t^{(1)}(\xi)$ and
$N_t^{(2)}(\xi)$ the integrals over the regions $|\tau| \leq
\sqrt{2}$, respectively $|\tau| \geq \sqrt{2}$. Since $|\sin x|
\leq x$ for all $x>0$, we obtain:
$$\frac{1}{a}|\cF_{0,t}\sin (a \cdot)(\tau)| \leq \int_{0}^{t}\frac{|\sin(as)|}{a}ds \leq \int_0^t s ds=\frac{t^2}{2},$$
and hence, by (\ref{Fourier-sin2}),
\begin{eqnarray}
\nonumber N_t^{(1)}(\xi)&=&\int_{|\tau| \leq \sqrt{2}}
\frac{1}{a^2}|\cF_{0,t} \sin (a \cdot)(\tau)|^2 \nu(d\tau)\leq \frac{t^4}{4}\int_{|\tau| \leq \sqrt{2}}\nu(d\tau) \\
\label{Nt-wave-UB1-step1} &\leq & \frac{3t^2}{4} \int_{|\tau| \leq
\sqrt{2}}\frac{1}{\tau^2+1}\nu(d\tau).
\end{eqnarray}

When $|\tau| \geq \sqrt{2}$, we have $\tau^2-a^2 \geq
\frac{1}{2}\tau^2$, since $a^2 \leq 1 \leq \frac{1}{2}\tau^2$.
Using the fact that $|f_{T}(\tau)| \leq 2T |\tau|$ and
$|g_{T}(\tau)| \leq 2$, we obtain:
$$f_{ta}^{2}\left(\frac{\tau}{a} \right)+g_{ta}^{2}\left(\frac{\tau}{a} \right) \leq 4 t^2 a^2 \frac{\tau^2}{a^2}+4=4(t^2 \tau^2+1).$$
Hence,
\begin{eqnarray}
\nonumber N_t^{(2)}(\xi) &=& \int_{|\tau| \geq
\sqrt{2}}\frac{1}{(\tau^2-a^2)^2}[f_{ta}^{2}\left(\frac{\tau}{a}
\right)+g_{ta}^{2}\left(\frac{\tau}{a} \right)]\nu(d\tau)
\leq  16 \int_{|\tau| \geq \sqrt{2}}\frac{t^2 \tau^2+1}{\tau^4} \nu(d\tau) \\
\nonumber
& \leq & 16\left(t^2+ \frac{1}{2}\right) \int_{|\tau| \geq \sqrt{2}}\frac{1}{\tau^2}\nu(d\tau) \\
\label{Nt-wave-UB1-step2}&\leq&  24\left(t^2+ \frac{1}{2}\right)
\int_{|\tau| \geq \sqrt{2}}\frac{1}{\tau^2+1}\nu(d\tau).
\end{eqnarray}
Relation (\ref{Nt-wave-UB1}) follows by taking the sum of
(\ref{Nt-wave-UB1-step1}) and (\ref{Nt-wave-UB1-step2}).

On the other hand, since $1 \leq b \leq \sqrt{2}$, $(|\tau|+b)^2
\leq (|\tau|+\sqrt{2})^2 \leq 4(\tau^2+1)$, and hence
\begin{equation}
\label{Nt-wave-UB2} \frac{1}{b} I(b) \geq \frac{1}{4\sqrt{2}}
\int_{\bR} \frac{1}{\tau^2+1}\nu(d\tau)=\frac{1}{4\sqrt{2}}K.
\end{equation}

From (\ref{Nt-wave-UB1}) and (\ref{Nt-wave-UB2}), we obtain:
\begin{equation}
\label{Nt-wave-UB5}
N_t(\xi) \leq 12 \sqrt{2} \max\{\frac{t^2}{4},8(t^2+\frac{1}{2}) \}\frac{1}{b}I(b)=:D_t \frac{I(b)}{b}.
\end{equation}

{\em (ii) Assume that $a \geq 1$.} In this case, $a \leq b \leq
\sqrt{2}a$. By Lemma \ref{hyp-UB-lemma1}, we have:
\begin{equation}
\label{Nt-wave-UB3} N_t(\xi) \leq C t \frac{1}{a}I(a) \leq C t
\frac{\sqrt{2}}{b}I(a).
\end{equation}

For any $x>0$, we write $I(x)=I_1(x)+I_2(x)$, where
$$I_1(x)=\int_{|\tau| \leq x} \frac{1}{(|\tau|+x)^2}\nu(d\tau) \quad \mbox{and} \quad
I_2(x)=\int_{|\tau| \geq x} \frac{1}{(|\tau|+x)^2}\nu(d\tau).$$
Note that
$$I_1(a) \leq  \frac{1}{a^2} \int_{|\tau| \leq a} \nu(d\tau) \leq \frac{2}{b^2} \int_{|\tau| \leq b}\nu(d\tau) \leq 8 I_1(b).$$
Using condition (C), we have:
\begin{eqnarray*}
I_2(a) & \leq & \int_{|\tau| \geq a}\frac{1}{\tau^2}\nu(d\tau) \leq \int_{|\tau| \geq b/\sqrt{2}} \frac{1}{\tau^2} \eta(|\tau|) d\tau=\sqrt{2} \int_{|\tau| \geq b}\frac{1}{\tau^2}\eta\left(\frac{1}{\sqrt{2}}|\tau|\right)d\tau \\
& \leq & \sqrt{2} C_{1/\sqrt{2}} \int_{|\tau| \geq b}
\frac{1}{\tau^2} \nu(d\tau) \leq 4 \sqrt{2} C_{1/\sqrt{2}} I_2(b).
\end{eqnarray*}
Therefore,
\begin{equation}
\label{Nt-wave-UB4} I(a) \leq \max\{8, 4 \sqrt{2}
C_{1/\sqrt{2}}\}I(b).
\end{equation}
From (\ref{Nt-wave-UB3}) and (\ref{Nt-wave-UB4}), we conclude
that:
\begin{equation}
\label{Nt-wave-UB6}
N_t(\xi) \leq C \sqrt{2} \max\{8, 4\sqrt{2}C_{1/\sqrt{2}}\}t \frac{1}{b}I(b)=:D t\frac{I(b)}{b}.
\end{equation}
The conclusion follows from (\ref{Nt-wave-UB5}) and (\ref{Nt-wave-UB6}) with $D_t^{(1)}=\max\{D_t,D t\}$.
$\Box$

We now treat the lower bound.

\begin{lemma}
\label{hyp-LB-lemma3} Assume that $\eta$ satisfies (C), and is
either non-increasing on $(0,\infty)$, or non-decreasing on
$(0,\infty)$. Then for any $t>0$ and $\xi \in \bR^d$, we have:
$$N_t(\xi) \geq C' t \frac{\eta(\sqrt{\Psi(\xi)})}{\Psi(\xi)}, \quad \mbox{if} \quad
t \sqrt{\Psi(\xi)} \geq 1,$$
where $C'$ is a positive constant.
\end{lemma}

\noindent {\bf Proof:} Let $a=\sqrt{\Psi(\xi)}$. Assume first that
$\eta$ is non-increasing on $(0,\infty)$.

Let $\varepsilon>0$ be arbitrary (to be chosen later). Since the
integrand of (\ref{def-Nt-hyp}) is non-negative, $N_t(\xi)$ is
bounded below by the integral over the region $\{|\tau| \leq
(1+\varepsilon)a\}$. In this region, using the fact that $\eta$ is
non-increasing and condition (C), we have
$$\eta(|\tau|) \geq \eta((1+\varepsilon)a) \geq C_{1/(1+\varepsilon)}^{-1} \eta(a)=:c_{\varepsilon}\eta(a).$$
Hence,
$$
N_t(\xi) \geq c_{\varepsilon} \eta(a)\int_{|\tau| \leq
(1+\varepsilon)a}\frac{1}{(\tau^2-a^2)^2}\left[f_{ta}^{2}
\left(\frac{\tau}{a} \right)+ g_{ta}^{2} \left(\frac{\tau}{a}
\right)\right]d\tau =: c_{\varepsilon}\eta(a)(I_1-I_2),$$ where
\begin{eqnarray*}
I_1&=& \int_{\bR} \frac{1}{(\tau^2-a^2)^2}\left[f_{ta}^{2} \left(\frac{\tau}{a} \right)+ g_{ta}^{2} \left(\frac{\tau}{a} \right)\right]d\tau\\
I_2&=& \int_{|\tau| \geq (1+\varepsilon)a}
\frac{1}{(\tau^2-a^2)^2}\left[f_{ta}^{2} \left(\frac{\tau}{a} \right)+
g_{ta}^{2} \left(\frac{\tau}{a} \right)\right]d\tau.
\end{eqnarray*}
Using Plancherel's theorem and the fact that $1-(\sin x)/x \geq
1/2$ for any $x \geq 2$, it follows that:
\begin{eqnarray*}
I_1&=&\frac{1}{a^2}\int_{\bR}|\cF_{0,t} \sin (a \ \cdot)(\tau)|^2 d\tau=\frac{1}{a^2} 2\pi \int_0^t \sin^2(as)ds \\
&=& \frac{1}{a^2} \pi t \left[ 1- \frac{\sin(2at)}{2at}\right]
\geq \frac{\pi}{2}\frac{t}{a^2}.
\end{eqnarray*}
On the other hand, using (\ref{important-ineq}), we obtain that:
$$I_2 \leq \frac{2t}{a} \int_{|\tau| \geq (1+\varepsilon)a}\frac{1}{(|\tau|-a)^2}d\tau=\frac{4t}{a} \int_{\varepsilon a}^{\infty}\frac{1}{\tau^2}d\tau=\frac{4}{\varepsilon}\frac{t}{a^2}.$$
Choose $\varepsilon>0$ such that
$\frac{\pi}{2}>\frac{4}{\varepsilon}$. Hence,
$$N_t(\xi) \geq c_{\varepsilon} \left(\frac{\pi}{2}-\frac{4}{\varepsilon}\right) t \frac{\eta(a)}{a^2}=:C' t \frac{\eta(a)}{a^2}.$$

If $\eta$ is non-decreasing on $(0,\infty)$, we bound $N_t(\xi)$
below by the integral over the region $|\tau| \leq
(1-\varepsilon)a$, for some $\varepsilon \in (0,1)$. In this
region,
$$\eta(|\tau|) \geq \eta((1-\varepsilon)a) \geq C_{1/(1-\varepsilon)}^{-1}\eta(a)=:c_{\varepsilon}'\eta(a).$$
As above, we obtain that
$N_t(\xi) \geq c_{\varepsilon}' \eta(a)(I_1-I_2')$,
where
\begin{eqnarray*}
I_2' &=& \int_{|\tau| \leq (1-\varepsilon)a}\frac{1}{(\tau^2-a^2)^2} \left[f_{ta}^{2} \left(\frac{\tau}{a} \right)+ g_{ta}^{2} \left(\frac{\tau}{a} \right)\right]d\tau \\
&=& \frac{2t}{a} \int_{|\tau| \leq
(1-\varepsilon)a}\frac{1}{(a-|\tau|)^2}d\tau=\frac{4t}{a}
\int_{\varepsilon
a}^{a}\frac{1}{\tau^2}d\tau=4\left(\frac{1}{\varepsilon}-1
\right)\frac{t}{a^2}.
\end{eqnarray*}
Choose $\varepsilon \in (0,1)$ such that
$\frac{\pi}{2}>4(\frac{1}{\varepsilon}-1)$, i.e.
$1>\varepsilon>(1+\frac{\pi}{8})^{-1}$. Hence,
$$N_t(\xi) \geq c_{\varepsilon} \left(\frac{\pi}{2}+4-\frac{4}{\varepsilon} \right)t \frac{\eta(a)}{a^2}=:C' t \frac{\eta(a)}{a^2}.$$
$\Box$

\begin{lemma}
\label{hyp-LB-lemma4} Suppose that $\eta$ satisfies either (C1) or
(C2). Then for any $K>0$ there exists a constant $M_K>0$ such that
$$\frac{1}{a} \int_{\bR}\frac{1}{(|\tau|+a)^2}\nu(d\tau) \leq M_K \frac{\eta(a)}{a^2}, \ \mbox{for any} \ a \geq K.$$
\end{lemma}

\noindent {\bf Proof:} If (C1) holds, then
\begin{eqnarray*}
& & \frac{1}{a} \int_{|\tau| \leq a}\frac{1}{(|\tau|+a)^2} \eta(|\tau|) d\tau \leq \frac{2}{a^3} \int_{0}^{a}\eta(\tau)d \tau \leq 2D_K \frac{\eta(a)}{a^2} \\
& & \frac{1}{a} \int_{|\tau| \geq a}\frac{1}{(|\tau|+a)^2}
\eta(|\tau|) d\tau \leq \frac{2 \eta(a)}{a^2}
\int_{a}^{\infty}\frac{1}{\tau^2}d \tau = 2 \frac{\eta(a)}{a^2}.
\end{eqnarray*}
If (C2) holds, then
\begin{eqnarray*}
& & \frac{1}{a} \int_{|\tau| \leq a}\frac{1}{(|\tau|+a)^2} \eta(|\tau|) d\tau \leq \frac{2 \eta(a)}{a} \int_{0}^{a} \frac{1}{a^2}d \tau = 2\frac{\eta(a)}{a^2} \\
& & \frac{1}{a} \int_{|\tau| \geq a}\frac{1}{(|\tau|+a)^2}
\eta(|\tau|) d\tau \leq \frac{2}{a}
\int_{a}^{\infty}\frac{1}{\tau^2}d \tau \leq 2 D_K
\frac{\eta(a)}{a^2}.
\end{eqnarray*}
$\Box$

\begin{lemma}
Assume that $\eta$ satisfies (C). Suppose in addition that $\eta$
satisfies either (C1) or (C2). Then, for any $t>0$ and $\xi \in
\bR^d$, we have:
$$N_t(\xi) \geq  D_t^{(2)} \frac{1}{\sqrt{1+\Psi(\xi)}}\int_{\bR} \frac{1}{(|\tau|+\sqrt{1+\Psi(\xi)})^2}\nu(d\tau),$$
where $D_t^{(2)}$ is a positive constant depending on $t$.
\end{lemma}

\noindent {\bf Proof:} We denote $a=\sqrt{\Psi(\xi)}$ and
$b=\sqrt{1+\Psi(\xi)}$.

{\em (i) Assume $at \leq 1$.} We use the fact that $\cos$ is
decreasing on $[0,\frac{\pi}{2}]$. Let $0<c<d<\frac{\pi}{2}$ be
such that $\cos c <\frac{1}{\sqrt{2}}\cos 1$.

Since the integrand of (\ref{def-Nt-hyp}) is non-negative, the
integral is bounded below by the integral over the region $\{c
\leq |\tau t| \leq d\}$. For any $\tau$ in this region,
$(\tau^2-a^2)^2 \leq 2(\tau^4+a^4) \leq 2(d^4+1)/t^4$ and
$$g_{ta}^2 \left(\frac{\tau}{a} \right)=[\cos(\tau t)-\cos(at)]^2 \geq \frac{1}{2} \cos^2(at)-\cos^2(\tau t) \geq \frac{1}{2}\cos^2 1-\cos^2 c=:M>0.$$
Hence,
\begin{equation}
\label{LB-a-small1} N_t(\xi) \geq  \int_{c \leq |\tau t| \leq
d}\frac{1}{(\tau^2-a^2)^2}g_{ta}^2 \left(\frac{\tau}{a}
\right)\nu(d\tau)
 \geq  M \frac{t^4}{2(t^4+1)} \int_{c/t \leq |\tau| \leq d/t}\nu(d\tau)=:A(t).
 \end{equation}

On the other hand, since $b \geq 1$,
\begin{equation}
\label{LB-a-small2}K=\int_{\bR}\frac{1}{\tau^2+1}\nu(d\tau) \geq
\int_{\bR}\frac{1}{(|\tau|+1)^2} \geq
\frac{1}{b}\int_{\bR}\frac{1}{(|\tau|+b)^2}\nu(d\tau).
\end{equation}
From (\ref{LB-a-small1}) and (\ref{LB-a-small2}), we obtain that:
$$N_t(\xi) \geq \frac{A(t)}{K} \frac{1}{b}\int_{\bR}\frac{1}{(|\tau|+b)^2}\nu(d\tau).$$

{\em (ii) Assume $at \geq 1$.} Using Lemma \ref{hyp-LB-lemma3} and
Lemma \ref{hyp-LB-lemma4} (with $K=1/t$), we get:
$$N_t(\xi) \geq C t \frac{\eta(a)}{a^2} \geq C t \frac{1}{M_{1/t}} \frac{1}{a}\int_{\bR} \frac{1}{(|\tau|+a)^2}\nu(d\tau) \geq
C t\frac{1 }{M_{1/t}} \frac{1}{b}\int_{\bR}
\frac{1}{(|\tau|+b)^2}\nu(d\tau).$$ $\Box$

This concludes the proof of the lower bound and the proof of
Theorem \ref{main-th-hyp}.

The next result is a by-product of the previous lemmas, and gives
an alternative condition for the existence of the random field
solution to (\ref{hyperbolic-eq}).

\begin{corollary}
\label{cor-hyp}
Under the conditions of Theorem \ref{main-th-hyp}, for any $t>0$
and $\xi \in \bR^d$,
$$K_t^{(2)} \frac{\eta(\sqrt{1+\Psi(\xi)})}{1+\Psi(\xi)} \leq N_t(\xi)\leq K_t^{(1)} \frac{\eta(\sqrt{1+\Psi(\xi)})}{1+\Psi(\xi)},$$
where $K_t^{(1)}$ and $K_t^{(2)}$ are positive constants depending
on $t$. Consequently, equation (\ref{hyperbolic-eq}) has a random
field solution if and only if
\begin{equation}
\label{cond-hyp} \int_{\bR^d}
\frac{\eta(\sqrt{1+\Psi(\xi)})}{1+\Psi(\xi)}\mu(d\xi)<\infty.
\end{equation}
\end{corollary}

\noindent {\bf Proof:} Let $a=\sqrt{\Psi(\xi)}$ and
$b=\sqrt{1+\Psi(\xi)}$. The upper bound follows from Lemma
\ref{hyp-UB-lemma2} and Lemma \ref{hyp-LB-lemma4} (applied to $b$).

For the lower bound, assume first that $at \leq 1$. Then $1 \leq b
\leq \sqrt{1+t^{-2}}$. Using (\ref{LB-a-small1}) and the monotonicity of $\eta$, we obtain that:
$$N_t(\xi)\geq A(t) \geq C(t)\frac{\eta(b)}{b^2},$$
where $C(t)$ is a constant depending on $t$.

Assume next that $at \geq 1$. Then $a \leq b \leq \sqrt{t^2+1} \
a$. By Lemma \ref{hyp-LB-lemma3} and using the fact that $\eta$ is
non-increasing (or $\eta$ is non-decreasing and satisfies (C)), we
obtain that:
$$N_t(\xi) \geq Ct \frac{\eta(a)}{a^2} \geq C'(t) \frac{\eta(b)}{b^2},$$
where $C'(t)$ is a constant depending on $t$. $\Box$

\begin{remark}
{\rm In the case of Examples \ref{Riesz-ex} and \ref{Bessel-ex}
with $\gamma \in (-1,1)$, condition (\ref{cond-hyp}) becomes:
\begin{equation}
\label{cond-hyp-Riesz} \int_{\bR^d}
\left(\frac{1}{1+\Psi(\xi)}\right)^{1+\gamma/2}\mu(d\xi)<\infty.
\end{equation}
When applied to Example \ref{Riesz-ex} with the parametrization
$\gamma=2H-1$, $H \in (0,1)$, Corollary \ref{cor-hyp} becomes an
extension of Theorem 3.1 of \cite{BT10} to the case $H<1/2$. Note
that if $\gamma>0$ (respectively $\gamma<0$), condition
(\ref{cond-hyp-Riesz}) is stronger (respectively weaker) than
(\ref{cond-par-Riesz}). }
\end{remark}

\appendix

\section{The stochastic integral with respect to $\tM$}
\label{app-stoch-integr}

Let $\tM=\{\tM(A);A \in \cR_d\}$ be a complex random measure on $\bR^d$ with orthogonal increments and control measure $\tmu$, as specified by Definition \ref{def-tM}.

Let $\varphi \in L_{\bC}^2(\bR^d,\tmu)$ be a function
which is continuous on compact sets. We give below the construction of the stochastic integral of $\varphi$ with respect to $\tM$.

\vspace{1mm}

{\em Step 1.} Let $A \in \cR_d$ be fixed. The stochastic integral of $\varphi$ over $A$, with respect to $\tM$ is defined as the $L_{\bC}^2(\Omega)$-limit of  $X_n=\sum_{j=1}^{k_n}\varphi(\tau_{j,n})\tM(A_{j,n})$, $n \geq 1$, where
$\Delta_n=(A_{j,n})_{1 \leq j \leq k_n}$ 
is a partition of $A$ into sets from $\cR_d$ such that $\|\Delta_n\|=\max_{j \leq k_n}|A_{j,n}| \to 0$, and $\tau_{j,n} \in A_{j,n}$ is arbitrary. The limit exists, since the sequence $(X_n)_{n \geq 1}$ is Cauchy: using the orthogonality and additivity of $\tM$, one can prove that
$E|X_n-X_m|^2=\int_{A}|\varphi_n-\varphi_m|^2 d\tmu$,
where $\varphi_n=\sum_{j=1}^{k_n}\varphi(\tau_{j,n})1_{A_{j,n}}$ and $\int_{A}|\varphi_n|^2 d\tmu \to \int_{A}|\varphi|^2 d\tmu$, by the definition of the Stieltjes integral (see e.g. p. 228 of \cite{billingsley95}). In this calculation, we used the fact that
$$E[\tM(A)\overline{\tM(B)}]=\tmu(A \cap B), \quad \mbox{for all} \ A,B \in \cR_d$$
which is again a consequence of the orthogonality and additivity of $\tM$.
One can show that the limit of $(X_n)_n$ does not depend on the choice of $(\Delta_n)_n$ and $(\tau_{j,n})_{j,n}$. This limit is denoted by $\tM_{A}(\varphi)=\int_{A}\varphi(\tau)\tM(d\tau)$.

This stochastic integral has the following properties:\\
(a) $\tM_{A}(\varphi+\psi)=\tM_A(\varphi)+\tM_{A}(\psi)$ a.s.;\\
(b) $E[\tM_A(\varphi)]=0$ and $E[\tM_{A}(\varphi)\overline{\tM_A(\psi)}]=\int_{A}\varphi \overline{\psi} d\tmu$; \\
(c) $\tM_{A}^{-}(\varphi)=\overline{\tM}_{A}(\varphi)$, where $\tM^{-}(A)=\tM(-A)$ and $\overline{\tM}(A)=\overline{\tM(A)}$.

\vspace{1mm}

{\em Step 2.} The stochastic integral of $\varphi$ with respect to $\tM$ is defined as the $L_{\bC}^2(\Omega)$-limit of $Y_n=\tM_{A_n}(\varphi)$, $n \geq 1$, where the sequence $(A_n)_n \subset \cR_d$ is chosen such that $A_n \subset A_{n+1}$ for all $n$, and $\cup_n A_n=\bR^d$. The sequence $(Y_n)_n$ is Cauchy, since $E|Y_n-Y_m|^2=\int_{A_m \verb2\2 A_n}|\varphi|^2 d\tmu$ for any $m>n$, and $\int_{A_n}|\varphi|^2 d\tmu \to \int_{\bR^d}|\varphi|^2 d\tmu$ by the monotone convergence theorem. Since the limit of $(Y_n)_n$ does not depend on the choice of $(A_n)_n$, we denote it by $\tM(\varphi)=\int_{\bR^d}\varphi(\tau)\tM(d\tau)$.
This stochastic integral enjoys properties similar to (a)-(c) above. In particular,
\begin{equation}
\label{cov-M}
E[\tM(\varphi)\overline{\tM(\psi)}]=\int_{\bR^{d}}\varphi\overline{\psi}d\tmu.
\end{equation}

Moreover, $\tM(\varphi)$ is a real-valued random variable, for any function $\varphi$ which satisfies $\overline{\varphi(\tau)}=\varphi(-\tau)$ for all $\tau \in \bR^d$. This follows since
\begin{equation}
\label{M-real-valued}\overline{\tM(\varphi)}=\int_{\bR^d}\overline{\varphi(\tau)}\overline{\tM}(d\tau)=\int_{\bR^d} \varphi(-\tau)\tM^{-}(d\tau)=\tM(\varphi),
\end{equation}
where we used from property (c) for the second equality above.

\section{Some elementary estimates}

\begin{lemma}
\label{lemmaA} Let $\gamma \in (-1,1)$ be arbitrary.
(i) For any $a>0$, we have:
$$\int_{\bR}\frac{1}{\tau^2+a^2}|\tau|^{-\gamma}d\tau =C_{\gamma} a^{-\gamma-1},$$
where $C_{\gamma}=\int_{\bR}(s^2+1)^{-1} s^{-\gamma}ds$.
(ii) For any $a>1$, we have:
$$C_{\gamma}^{(1)} a^{-\gamma-1} \leq \int_{\bR}\frac{1}{\tau^2+a^2}(\tau^2+1)^{-\gamma/2}d\tau \leq C_{\gamma}^{(2)} a^{-\gamma-1},$$
where $C_{\gamma}^{(1)}$ and $C_{\gamma}^{(2)}$ are some positive
constants depending on $\gamma$.
\end{lemma}

\noindent {\bf Proof:} (i) The conclusion follows by the change of
variable $\tau/a=\tau'$.

(ii) We denote by $I_1$ and $I_2$ the integrals over the regions
$|\tau| \leq a$, respectively $|\tau| \geq a$. We use the notation
$f(\tau) \sim g(\tau)$ if $c_1 g(\tau)\leq f(\tau) \leq c_2
g(\tau)$ for any $\tau \in \bR$, for some constants $c_1,c_2>0$.
When $|\tau| \leq a$,
$$(\tau^2+a^2)^{-1} (\tau^2+1)^{-\gamma/2} \sim a^{-2} (\tau^2+1)^{-\gamma/2} \sim
a^{-2} (\tau+1)^{-\gamma},$$ and hence $I_1$ is bounded below and
above by some constants multiplied by:
$$a^{-2}\int_{0}^{a}(\tau+1)^{-\gamma}d\tau=a^{-2} \frac{1}{1-\gamma}[(a+1)^{1-\gamma}-1] \sim a^{-2}(a+1)^{1-\gamma} \sim a^{-\gamma-1}.$$

When $|\tau| \geq a$,
$$(\tau^2+a^2)^{-1} (\tau^2+1)^{-\gamma/2} \sim \tau^{-2} (\tau^2+1)^{-\gamma/2} \sim
\tau^{-2} (\tau+1)^{-\gamma} \sim \tau^{-2-\gamma},$$ and hence
$I_2$ is bounded below and above by some constants multiplied by:
$$\int_{a}^{\infty}\tau^{-2-\gamma}d\tau=\frac{1}{\gamma+1}a^{-\gamma-1}.$$
$\Box$

\section{Verification of conditions (hyperbolic case)}

\begin{lemma}
\label{lemma-B1}
(i) Let $\eta(\tau)=|\tau|^{-\gamma}$ with $\gamma \in (-1,1)$.
Then $\eta$ satisfies (C). If $\gamma \in (0,1)$, then $\eta$
satisfies (C1). If $\gamma \in (-1,0)$, then $\eta$ satisfies
(C2).

(ii) Let $\eta(\tau)=(1+\tau^2)^{-\gamma/2}$ with $\gamma >-1$. Then
$\eta$ satisfies (C). If $\gamma \in (0,1)$, then $\eta$ satisfies
(C1). If $\gamma \in (-1,0)$, then $\eta$ satisfies (C2).
\end{lemma}

\noindent {\bf Proof:} (i) (C) is clearly satisfied since
$\eta(\lambda \tau)/\eta(\tau)=\lambda^{-\gamma}=:C_{\lambda}$ for
any $\tau>0$.

If $\gamma \in (0,1)$, (C1) holds since for any $a>0$,
$$\int_0^a \eta(\tau)d\tau=\int_0^a \tau^{-\gamma}d\tau=\frac{1}{1-\gamma}a^{1-\gamma}=\frac{1}{1-\gamma}a \eta(a).$$

If $\gamma \in (-1,0)$, (C2) holds since for any $a>0$,
$$\int_{a}^{\infty} \tau^{-2}
\eta(\tau)d\tau=\int_{a}^{\infty}\tau^{-2-\gamma}d\tau=\frac{1}{\gamma+1}a^{-1-\gamma}=
\frac{1}{\gamma+1} a^{-1}\eta(a).$$

(ii) We first check (C). If $\gamma>0$, then the inequality
\begin{equation}
\label{ineq-check-C} \frac{\eta(\lambda
\tau)}{\eta(\tau)}=\left(\frac{1+\tau^2}{1+\lambda^2
\tau^2}\right)^{\gamma/2} \leq C_{\lambda}
\end{equation}
is equivalent to $1+\lambda^2 \tau^2 \geq
C_{\lambda}^{-2/\gamma}(1+\tau^2)$. We choose
$C_{\lambda}=\{\min(1,\lambda^2)\}^{-\gamma/2}$.

If $\gamma<0$, then (\ref{ineq-check-C}) is equivalent to
$1+\lambda^2 \tau^2 \leq C_{\lambda}^{-2/\gamma}(1+\tau^2)$. We
choose $C_{\lambda}=\{\max(1,\lambda^2)\}^{-\gamma/2}$. If $\gamma
\in (0,1)$ then (C1) holds since for any $a \geq K$,
$$\int_0^a (1+\tau^2)^{-\gamma/2}d\tau \leq \int_0^a \tau^{-\gamma}d\tau=\frac{a}{1-\gamma}\left(\frac{1}{a^2} \right)^{\gamma/2} \leq \frac{a}{1-\gamma} \left(\frac{C_K}{1+a^2}\right)^{\gamma/2},$$
where $C_K=1+K^{-2}$. If $\gamma \in (-1,0)$ then (C2) holds since
for any $a \geq K$,
\begin{eqnarray*}
& & \int_a^{\infty} \tau^{-2}(1+\tau^2)^{-\gamma/2}d\tau \leq
\int_{a}^{\infty}\tau^{-2}
(C_K \tau^2)^{-\gamma/2}=C_{K}^{-\gamma/2}\frac{1}{\gamma+1}a^{-\gamma-1} \\
& & \ \ \ \ \ \ =C_{K}^{-\gamma/2}\frac{1}{\gamma+1}a^{-1} \left(
\frac{1}{a^2}\right)^{\gamma/2} \leq C_{K}^{-\gamma/2}a^{-1}
\left(\frac{1}{1+a^2}\right)^{\gamma/2}.
\end{eqnarray*}
$\Box$

\vspace{3mm}

\footnotesize{{\em Acknowledgement.} The author would like
to thank Davar Khoshnevisan for drawing her attention to the class of processes with stationary increments as a possible replacement for the fBm (which appears in the temporal component of the Gaussian noise of \cite{B11}),
and also for pointing out a simplified proof of Lemma \ref{lemmaA}.(i).

\normalsize{

\end{document}
\begin{thebibliography}{99}

\bibitem{B11} Balan, R. M. (2011). Some linear SPDEs driven by a fractional noise with Hurst index greater than $1/2$. Preprint available on arXiv:1102.3992.

\bibitem{BT08}  Balan, R.M. and Tudor, C. A. (2008).
The stochastic heat equation with fractional-colored noise:
existence of the solution. {\em Latin Amer. J. Probab. Math.
Stat.} {\bf 4}, 57-87.



\bibitem{BT10} Balan, R. M. and Tudor, C. A. (2010). The stochastic
wave equation with fractional noise: a random field approach. {\em
Stoch. Proc. Appl.} {\bf 120}, 2468-2494.

\bibitem{bardina-florit05} Bardina, X. and Florit, C. (2005). Approximation in law to the $d$-parameter FBS based on the functional invariance principle. {\em Rev. Mat. Iberoamer.}
{\bf 21}, 1037-1052.

\bibitem{billingsley95} Billingsley, P. (1995).
{\em Probability and Measure}. Third Edition. John Wiley, New
York.

\bibitem{bonami-estrade03} Bonami, A. and Estrade, A. (2003). Anisotropic analysis of some
Gaussian models. {\em J. Fourier Anal. Appl.} {\bf 9}, 215-236.

\bibitem{clausel-vesel11} Clausel, M. and Vesel, B. (2011). Explicit construction of
operator scaling Gaussian random fields. Preprint available on arXiv:1104.0774.


\bibitem{dalang99} Dalang, R. C. (1999). Extending martingale
measure stochastic integral with applications to spatially
homogenous s.p.d.e.'s. {\em Electr. J. Probab.} {\bf 4}, no. 6, 29
pp.

\bibitem{dalang-mueller03} Dalang, R. C. and Mueller, C. (2003). Some non-linear S.P.D.E.'s
that are second order in time. {\em Electr. J. Probab.} {\bf 8},
no. 1, 1-21.

\bibitem{dalang-sanzsole05} Dalang, R. C. and Sanz-Sol\'e, M. (2005). Regularity of the sample paths of a class of second order spde's. {\em J. Funct. Anal.} {\bf 227}, 304-337.

\bibitem{dobrushin-major79} Dobrushin, R. L. and Major, P. (1979).
Non-central limit theorems for non-linear functionals of
 Gaussian fields. {\em Z. Wahrsch. verw. Gebiete}
{\bf 50}, 27-52.

\bibitem{doob53} Doob, J. L. (1953). {\em Stochastic Processes}.
John Wiley, New York.

\bibitem{FK10} Foondun, M. and Khoshnevisan, D. (2010). On the stochastic heat equation with spatially-colored random forcing. Preprint available on arXiv:1003.0348.

\bibitem{hu01} Hu, Y. (2001). Heat equations with fractional white
noise potentials. {\em Appl. Math. Optim.} {\bf 43}, 221-243.

\bibitem{hu-nualart-song11} Hu, Y., Nualart, D. and Song, J. (2011). Feynman-Kac formula for heat equation driven by fractional white noise. {\em Ann. Probab.} {\bf 39}, 291-326.


\bibitem{jolis07} Jolis, M. (2007). On the Wiener integral with
respect to the fractional Bronwian motion on an interval. {\em J.
Math. Anal. Appl.} {\bf 330}, 1115-1127.

\bibitem{jolis10} Jolis, M. (2010). The Wiener integral with
respect to second order processes with stationary increments. {\em
J. Math. Anal. Appl.} {\bf 366}, 607-620.


\bibitem{kamont96} Kamont, A. (1996). On the fractional anisotropic Wiener field. {\em Probab. Math. Stat.} {\bf 18}, 85-98.

\bibitem{K02} Khoshnevisan, D. (2002). {\em Multiparameter Processes. An Introduction to Random Fields}. Springer-Verlag, New York.

\bibitem{KX09} Khoshnevisan, D. and Xiao, Y. (2009). Harmonic analysis of additive L\'evy processes. {\em Probab. Th. Rel. Fields} {\bf 145}, 459-515.

\bibitem{leger-pontier99} L\'eger, S. and Pontier, M. (1999). Drap brownien fractionnaire. {\em CRAS Paris Serie I}. {\bf 329}, 893-898.

\bibitem{li-xiao11} Li, Y. and Xiao, Y. (2011). Multivariate operator-self-similar random fields. Preprint available on arXiv:1104.0059.


\bibitem{lindstrom93} Lindstrom, M. (1993). Fractional Brownian fields as integrals of
white noise. {\em Bull. London Math. Soc.} {\bf 25}, 83-88.

\bibitem{mura-mainardi09} Mura, A. and Mainardi, F. (2009). A class of self-similar stochastic processes with stationary increments to model anomalous diffusion in physics. To appear in {\em Integr. Transforms Special Funct.} {\bf 20}. Proceedings of the conference ``Linear and non-linear theory of generalized functions and its applications'', Bedlewo, Poland, 2007.

\bibitem{PT00} Pipiras, V. and Taqqu, M. S. (2000). Integration questions
related to fractional Brownian motion. {\em Probab. Th. Rel.
Fields} {\bf 118}, 251-291.


\bibitem{PT01} Pipiras, V. and Taqqu, M. (2001).
Are classes of deterimistic integrands for the fractional Brownian
motion on a finite interval complete? {\em Bernoulli} {\bf 7},
873--897.

\bibitem{stein70} Stein, E. M. (1970). {\em Singular Integrals and Differentiability Properties of Functions}. Princeton University Press. Princeton, New Jersey.

\bibitem{yaglom57} Yaglom, A. M. (1957). Some classes of random
fields in $n$-dimensional space, related to stationary random
processes. {\em Th. Probab. Appl.} {\bf 2}, 273-320.


\bibitem{yaglom87} Yaglom, A. M. (1987).
{\em Correlation Theory of Stationary and Related Random
Functions}. Springer.


\end{thebibliography}
